\documentclass[a4paper,12pt, reqno]{amsart}

\usepackage{amsmath}
\usepackage{amssymb}
\usepackage{mathrsfs}
\usepackage{ifthen}
\usepackage{graphicx}
\usepackage[T1]{fontenc} 
\usepackage{lipsum}
\usepackage[a4paper,margin=0.75in]{geometry}

\usepackage{cite}

\def\switchlinenumbers{\@ifstar
	{\let\makeLineNumberOdd\makeLineNumberRight
		\let\makeLineNumberEven\makeLineNumberLeft}%
	{\let\makeLineNumberOdd\makeLineNumberLeft
		\let\makeLineNumberEven\makeLineNumberRight}%
}

\def\setmakelinenumbers#1{\@ifstar
	{\let\makeLineNumberRunning#1%
		\let\makeLineNumberOdd#1%
		\let\makeLineNumberEven#1}%
	{\ifx\c@linenumber\c@runninglinenumber
		\let\makeLineNumberRunning#1%
		\else
		\let\makeLineNumberOdd#1%
		\let\makeLineNumberEven#1%
		\fi}%
}

\nonstopmode \numberwithin{equation}{section}
\newtheorem*{theorem*}{Theorem}

\newtheorem{thm}{Theorem}[section]
\newtheorem{cor}{Corollary}[section]
\newtheorem{lem}{Lemma}[section]

\theoremstyle{definition}
\newtheorem{defn}{Definition}[section]
\newtheorem{example}{Example}[section]
\newtheorem{qsn}{Question} [section]
\newtheorem{prob}[equation]{Problem}
\newtheorem{rem}{Remark}[section]

\newtheorem{innercustomthm}{Theorem}
\newenvironment{customthm}[1]
{\renewcommand\theinnercustomthm{#1}\innercustomthm}
{\endinnercustomthm}

\newcounter{minutes}
\divide\time by 60
\newcounter{hours}
\multiply\time by 60
\addtocounter{minutes}{-\time}

\newcounter {own}
\def\theown {\thesection       .\arabic{own}}

\newenvironment{pf}[1][]{%
	\vskip 3mm
	\noindent
	\ifthenelse{\equal{#1}{}}%
	{{\slshape {\bf Proof}. }}%
	{{\slshape #1.} }%
}%
{\qed\bigskip}

\newcounter{alphabet}

\def\be{\begin{equation}}
	\def\ee{\end{equation}}

\newcommand{\bee}{\begin{enumerate}}
	\newcommand{\eee}{\end{enumerate}}

\newcommand{\blem}{\begin{lem}}
	\newcommand{\elem}{\end{lem}}
\newcommand{\bthm}{\begin{thm}}
	\newcommand{\ethm}{\end{thm}}
\newcommand{\bcor}{\begin{cor}}
	\newcommand{\ecor}{\end{cor}}
\newcommand{\beg}{\begin{examp}}
	\newcommand{\eeg}{\end{examp}}
\newcommand{\begs}{\begin{examples}}
	\newcommand{\eegs}{\end{examples}}
\newcommand{\bdefe}{\begin{defin}}
	\newcommand{\edefe}{\end{defin}}
\newcommand{\bprob}{\begin{prob}}
	\newcommand{\eprob}{\end{prob}}
\newcommand{\bei}{\begin{itemize}}
	\newcommand{\eei}{\end{itemize}}
\newcommand{\real}{{\operatorname{Re}\,}}

\newcommand{\norm}[1]{\left\lVert#1\right\rVert}

\newcommand{\comment}[1]{}

\subjclass[{AMS} Subject Classification:]{Primary 32A05, 31C10, 46B07;  Secondary 32Q02, 46E40}

\begin{document}
	
	\title[]{Local Banach Space Theoretic Approach to Bohr's Theorem for Vector Valued Holomorphic and Pluriharmonic Functions}

	\author{Himadri Halder}
	\address{Himadri Halder,
			Statistics and Mathematics Unit,
			Indian Statistical Institute, Bangalore Centre,
			8th Mile Mysore Road, R. V. College Post, Bangalore-560059, India}
	\email{himadrihalder119@gmail.com}
	

	
	
	\begin{abstract}
		We study Bohr's theorem for vector valued holomorphic and operator valued pluriharmonic functions on complete Reinhardt domains in $\mathbb{C}^n$. Using invariants from local Banach space theory, we show that the associated Bohr radius is always strictly positive and obtain its asymptotic behavior separately in the finite- and infinite-dimensional settings. 
		Our analysis provides asymptotic estimates for this constant on convex as well as non-convex complete Reinhardt domains.
		The framework developed here includes the classical Minkowski-space setting as a special case and applies to a wide class of Banach sequence spaces, including mixed Minkowski, Lorentz, and Orlicz spaces. We further establish a coefficient-type Schwarz–Pick lemma for operator valued pluriharmonic maps on complete Reinhardt domains.	
	\end{abstract}

	\maketitle
	\pagestyle{myheadings}
	\markboth{}{A Local Banach Space Theoretic Approach to Bohr’s Theorem}
	
	\section{Introduction and the main results}\label{section-1}
	This paper aims to contribute to the study of Bohr's theorem for vector valued holomorphic and operator valued pluriharmonic functions through the framework of local Banach space theory, an area not previously explored. 
	The framework and unified approach developed here provide a powerful tool for studying both vector valued holomorphic and pluriharmonic functions in broad classes of Banach sequence spaces, representing one of the most significant applications of our results in functional and complex analysis. 
	Over the last three decades, Bohr's theorem and its applications have occupied a central position in geometric function theory, connecting ideas from several complex variables, functional analysis, and operator theory. For a given complete Reinhardt domain $\Omega \subset \mathbb{C}^n$, the Bohr radius $K_{n}(\Omega)$ of $\Omega$ is defined to be the supremum of all $r \in [0,1]$ such that for each holomorphic function $f(z)=\sum_{\alpha}c_{\alpha}\,z^{\alpha}:\Omega\rightarrow \mathbb{C}$ we have
	\begin{equation} \label{e-hh-mult-BI}
		\sup_{z \in r\Omega}\, \sum_{\alpha}|c_{\alpha}\,z^{\alpha}| \leq \sup_{z \in \Omega}\, \left|f(z)\right|.
	\end{equation}
	Note that with this notation, Bohr's celebrated power series theorem states that $K_{1}(\mathbb{D})=1/3$ (see \cite{Bohr-1914}). One of the most striking aspects of Bohr's theorem arises in multidimensional settings, where the situation becomes significantly more complex and many problems remain open. In particular, there is still limited knowledge about the behavior of the constant $K_{n}(\Omega)$. A key challenge lies in determining the constant $K_{n}(\Omega)$, whose exact value remains unknown for $n>1$, even in the polydisc $\mathbb{D}^n$, where $\mathbb{D}^n:=\{z=(z_{1},\ldots,z_{n}) \in \mathbb{C}^n: |z_{j}|<1\, \mbox{for all}\, 1 \leq j \leq n\}$. To advance this direction, researchers have focused on bounding $K_{n}(\Omega)$ and analyzing its asymptotic behavior on general domains, and in particular on key domains such as the unit ball of Minkowski space and finite-dimensional complex Banach spaces. Beginning with the work of Dineen and Timoney \cite{Dineen-Timoney-1989}, further developed by Boas and Khavinson \cite{boas-1997,boas-2000}, and later extended by Aizenberg \cite{aizn-2000a} and by Defant and Frerick \cite{defant-2011-lpn}, it is known that for every $1\leq q \leq \infty$ and all $n\in \mathbb{N}$, there exist constants $C,D>0$ such that 
	\begin{equation*}
		\frac{1}{C}\, \left(\frac{\log\,n}{n}\right)^{1-\frac{1}{\min\{q,2\}}} \leq K_{n}(B_{\ell^n _q}) \leq D \,\left(\frac{\log\,n}{n}\right)^{1-\frac{1}{\min\{q,2\}}},
	\end{equation*}
	where $B_{\ell^n _q}:=\left\{z=(z_{1},\ldots,z_{n}) \in \mathbb{C}^n:\norm{z}_{q}:=\left(\sum_{i=1}^{n}|z_{i}|^q\right)^{1/q}<1\right\}$, $1\leq q <\infty$ and $B_{\ell^n _\infty}:=\mathbb{D}^n$.
	\vspace{2mm}
	
	A class of holomorphic functions $f$, defined in a complete Reinhardt domain $\Omega$, is said to exhibit the Bohr phenomenon on $\Omega$ if there exists a universal constant $r=r_{0}\in (0,1]$, called the Bohr radius of $\Omega$ with respect to that class, such that the inequality \eqref{e-hh-mult-BI} holds for all functions in that class. Since not every class of functions exhibits the Bohr phenomenon (see \cite{aizn-2000b,bene-2004,Blasco-Collect-2017}), it is of interest to determine when a given class does. One of the main aims of this note is to address this question for vector valued functions within the framework of local Banach space theory, by establishing conditions under which the Bohr phenomenon can be verified in this setting. With this goal in mind, we first recall key facts about Banach space valued holomorphic functions and, for this purpose, introduce the Bohr radius for this class. Let $\Omega\subset \mathbb{C}^n$ be a complete Reinhardt domain and $n\in \mathbb{N}$. Let $U:X\rightarrow Y$ be a bounded liner operator between two complex Banach spaces and $\norm{U} \leq \lambda$. For $1 \leq p < \infty$, the $\lambda$-powered Bohr radius of $U$, denoted by $K_{\lambda}(\Omega, p,U)$, is defined to be the supremum of all $r\geq 0$ such that for all holomorphic functions $f(z)=\sum_{\alpha}a_{\alpha}z^{\alpha}:\Omega \rightarrow X$ we have 
	\begin{equation} \label{e-1.2}
		\sup_{z \in r\Omega}\, \sum_{\alpha} \norm{U(a_{\alpha})z^\alpha}^p_{Y} \leq \lambda^p\,\norm{f}^p_{\Omega,X}, 
	\end{equation}
	where $\norm{f}_{\Omega,X}:=\sup_{z \in \Omega}\norm{f(z)}_{X}$. Set $K(\Omega, p,U):=K_{1}(\Omega, p,U)$, $K_{\lambda}(\Omega, p,X):=K_{\lambda}(\Omega, p,U)$ whenever $U=I:X\rightarrow X$, $K(\Omega, p,X):=K_{1}(\Omega, p,X)$, $K_{\lambda}(\Omega, p):=K_{\lambda}(\Omega, p,\mathbb{C})$, and $K(\Omega,p):=K_1(\Omega,p)$. Note that with this notations, $K(\Omega, 1)$ coincides with $K_{n}(\Omega)$, Bohr radius for scalar valued case defined by \eqref{e-hh-mult-BI}. The above definition is inspired by Defant {\it et al.} \cite{defant-2012}, extending their case $p=1$ and $\Omega=\mathbb{D}^n$ to a broader setting, thereby enabling the study of the Bohr phenomenon for more general domains and function classes. Bohr's theorem for operator valued holomorphic functions defined in the unit disc was first investigated by Paulsen and Singh \cite{paulsen-2006}.
	\vspace{2mm}
	
	The constant $K_{\lambda}(\Omega, 1)$ was first studied by Bombieri \cite{bombieri-1962}, who determined its exact value for $\lambda \in [1,\sqrt{2}]$. Later, Bombieri and Bourgain \cite{bombieri-2004} have investigated its exact asymptotic behavior as $\lambda\rightarrow \infty$. Following these results, Defant {\it et al.} have introduced the more general constant $K_\lambda(\mathbb{D}^n,1,U)$ in \cite{defant-2012} and established its asymptotic estimates in both finite and infinite dimensional Banach spaces $X$. On a related note, the constants $K(\mathbb{D},p)$ and $K(\mathbb{D}^n,p)$ were first examined by Djakov and Ramanujan \cite{Djakov & Ramanujan & J. Anal & 2000}, and later developed further by B\'{e}n\'{e}teau, Dahlner, Khavinson, Das \cite{bene-2004}. In the Banach space valued setting, Blasco \cite{Blasco-Collect-2017} has shown that $K(\mathbb{D},p,X)=0$ for $p \in [1,2)$, while $K(\mathbb{D},p,X)>0$ if and only if $X$ is $p$-uniformly $\mathbb{C}$-convex for $2 \geq p < \infty$. 
	These results demonstrate that $K(\mathbb{D},p,X)$ is not necessarily nonzero for all Banach spaces $X$. This limitation motivated a refinement of its definition to ensure non-vanishing behavior for all $X$, which led to the introduction of the parameter $\lambda$ in \cite[Definition 1.1]{defant-2012}. With these developments in mind, it is natural to ask the following question.
	\begin{qsn} \label{qsn-1.1}
		For a given complete Reinhardt domain $\Omega$, can the the constant $K_{\lambda}(\Omega, p,U)$ be studied in detail? If so, does it always remain nonzero, and what can be said about its asymptotic behavior?
	\end{qsn}
	We answer affirmatively to this question in this paper. We now focus on complex-valued and operator-valued pluriharmonic functions defined on domains $\Omega$ in $\mathbb{C}^n$. A twice continuously differentiable function $f:\Omega \rightarrow \mathbb{C}$ is said to be plurharmonic if its restriction to every complex line is harmonic. Moreover, if $\Omega$ is a simply connected domain in $\mathbb{C}^n$ containing the origin, then $f:\Omega \rightarrow \mathbb{C}$ is pluriharmonic if it admits a representation $f=h+\overline{g}$, where $h,g$ are holomorphic on $\Omega$ with $g(0)=0$. Similarly, for a simply connected complete Reinhardt domain $\Omega$, a continuous function $f:\Omega \rightarrow \mathcal{B}(\mathcal{H})$ is pluriharmonic if and only if $f$ can be expressed as
	\begin{equation} \label{e-1.3-a}
		f(z)=\sum_{m=0}^{\infty} \sum_{|\alpha|=m} a_{\alpha}\, z^{\alpha} + \sum_{m=1}^{\infty} \sum_{|\alpha|=m} b^{*}_{\alpha}\, \bar{z}^{\alpha},
	\end{equation} 
	where $h(z)=\sum_{m=0}^{\infty} \sum_{\alpha} a_{\alpha}\, z^{\alpha}$ and $g(z)=\sum_{m=1}^{\infty} \sum_{\alpha} b_{\alpha}\, z^{\alpha}$ are $\mathcal{B}(\mathcal{H})$ valued holomorphic functions on $\Omega$. Here $\mathcal{B}(\mathcal{H})$ is the algebra of all bounded linear operators on a complex Hilbert space $\mathcal{H}$. For a simply connected complete Reinhardt domain, let $\mathcal{PH}(\Omega,X)$ denote the set of all bounded $X$ valued pluriharmonic functions on $\Omega$, where $X=\mathcal{B}(\mathcal{H})$. Here bounded means $\norm{f}_{\Omega,X}:=\sup_{z \in \Omega}\,\norm{f(z)}_{X} < \infty$ for $f \in \mathcal{PH}(\Omega,X)$.
	The Bohr radius problem for complex valued harmonic functions was first explored in \cite{Abu-2010}, later extended to operator valued functions in \cite{bhowmik-2021}, and more recently to pluriharmonic Hilbert space valued functions in \cite{hamada-Math-Nachr-2023}. These advances naturally lead to the following question for operator-valued pluriharmonic functions.
	\begin{qsn} \label{qsn-1.2}
		Let $\Omega \subset \mathbb{C}^n$ be a simply connected complete Reinhardt domain. Can we study the Bohr radius problem for operator valued pluriharmonic functions analogously to the inequality \eqref{e-1.2} and the constant $K_{\lambda}(\Omega, p,U)$?
	\end{qsn}
	In this paper, we provide an affirmative answer to this question, and to this end, we first introduce the notion of the powered Bohr radius for operator-valued pluriharmonic functions. 
	\begin{defn} \label{def-1.1}
		Let $X=\mathcal{B}(\mathcal{H})$ and $Y$ be any complex Banach space. Let $\Omega\subset \mathbb{C}^n$ be a simply connected complete Reinhardt domain and $n\in \mathbb{N}$. Let $U:X\rightarrow Y$ be a bounded liner operator and $\norm{U} \leq \lambda$. For $1 \leq p < \infty$, the $\lambda$-powered Bohr radius of $U$, denoted by $R_{\lambda}(\Omega, p,U)$, is defined to be the supremum of all $r\geq 0$ such that for all $f \in \mathcal{PH}(\Omega,X)$ of the form \eqref{e-1.3-a} we have 
		\begin{equation} \label{e-1.4-a}
			\sup_{z \in r\Omega}\,\sum_{m=0}^{\infty} \sum_{|\alpha|=m} (\norm{U(a_{\alpha})}^p_{Y} + \norm{U(b_{\alpha})}^p_{Y})|z^\alpha|^p \leq \lambda^p\,\norm{f}^p_{\Omega,X}. 
		\end{equation}
	\end{defn}
	\noindent Set $R(\Omega, p,U):=R_{1}(\Omega, p,U)$, $R_{\lambda}(\Omega, p,X):=R_{\lambda}(\Omega, p,U)$ whenever $U=I:X\rightarrow X$, $R(\Omega, p,X):=R_{1}(\Omega, p,X)$, $R_{\lambda}(\Omega, p):=R_{\lambda}(\Omega, p,\mathbb{C})$, and $R(\Omega,p):=R_1(\Omega,p)$. The asymptotic behaviour of the constants $R(\Omega,1)$ and $R(\Omega,p)$ have been studied in \cite{hamada-JFA-2022,das-2024} whenever $\Omega=B_{\ell^n _q}$.
	
	We now proceed to present our main results. To the best of our knowledge, no attempt has yet been made to study operator-valued analogues of Bohr’s theorem for pluriharmonic functions in the framework of Definition \ref{def-1.1}, not even in the simplest cases such as $U=I:X\rightarrow X$, $p=1$, $\lambda=1$, and for the classical domains such as the unit balls of Minkowski spaces. This gap motivates us to pose the following question.
	\begin{qsn} \label{qsn-1.3}
		For a given simply connected complete Reinhardt domain $\Omega$, can the the constant $R_{\lambda}(\Omega, p,U)$ be studied in detail? If so, does it always remain nonzero, and what can be said about its asymptotic behavior?
	\end{qsn} 
	In this paper, we address this question and begin by examining the constant $R_{\lambda}(\Omega, p,U)$ through invariants from local Banach space theory, such as unconditional bases and projection constants. Let $Z=(\mathbb{C}^n, ||.||)$ be an $n$-dimensional Banach space for which the canonical basis vectors $e_{k}$ form a normalized $1$-unconditional basis (or equivalently, the open unit ball $B_{Z}$ is a complete Reinhardt domain). There is a one-to-one correspondence between bounded convex complete Reinhardt domains $\Omega \subseteq\mathbb{C}^n$ and the open unit balls of norms on $\mathbb{C}^n$ for which the canonical basis vectors $e_{k}$ form a normalized $1$-unconditional basis. Indeed, the Minkowski functional $p_{\Omega}:\mathbb{C}^n \rightarrow \mathbb{R}_{+}$ of $\Omega$ defined by $$p_{\Omega}(z):= \inf \left\{t>0 : z/t \in \Omega\right\}$$
	 is a norm on $\mathbb{C}^n$, and $\Omega$ coincides with the open unit ball of $(\mathbb{C}^n,p_{\Omega})$. Conversely, if $Z=(\mathbb{C}^n, ||.||)$ is a Banach space whose canonical basis is normalized $1$-unconditional, then its open unit ball $B_{	Z}$ is clearly a bounded convex complete Reinhardt domain in $\mathbb{C}^n$. Thus, the study of the constant $R_{\lambda}(\Omega, p,U)$ for the unit balls $\Omega=B_{Z}$ of finite dimensional complex Banach spaces $(\mathbb{C}^n, ||.||)$ with the normalized $1$-unconditional canonical bases, is equivalent to studying $R_{\lambda}(\Omega, p,U)$ over bounded complete Reinhardt domains $\Omega \subseteq\mathbb{C}^n$.
	In the following result, we first establish that $R_{\lambda}(\Omega, p,U)$ is strictly positive whenever $\Omega=B_{Z}$. Combining this with the above correspondence and Lemma \ref{lem-3.5}, we conclude that $R_{\lambda}(\Omega, p,U)$ is strictly positive for every bounded simply connected complete Reinhardt domain $\Omega \subseteq \mathbb{C}^n$, even in the non-convex case.
	\begin{thm} \label{thm-1.1}
		Let $X=\mathcal{B}(\mathcal{H})$ and $U: X \rightarrow Y$ be a non-null bounded linear operator such that $\norm{U}< \lambda$. Then, for $\lambda>1$ and $n \in \mathbb{N}$, we have
		$R_{\lambda}(B_{Z}, p,U) \geq D. \frac{1}{\sup_{z \in B_{Z}}\norm{z}_{p}}$, where 
		$$
		D=\begin{cases}
			\max \left\{\frac{1}{4 \lambda\,2^{\frac{1}{p}}}\left(\frac{\lambda^p - \norm{U}^p}{2\lambda^p - \norm{U}}\right)^{\frac{1}{p}}\, , \frac{1}{4 \lambda\,2^{\frac{1}{p}}}\left(\frac{\lambda^p - \norm{U}^p}{\lambda^p - \norm{U}^p +1}\right)^{\frac{1}{p}}\, \frac{1}{\norm{U}}\right\}\,  & \text{for $\norm{U}\geq \frac{1}{4 \lambda\,2^{\frac{1}{p}}}$},\\[3mm]
			\max \left\{\frac{1}{4 \lambda\,2^{\frac{1}{p}}}\left(\frac{\lambda^p - \norm{U}^p}{2\lambda^p - \norm{U}}\right)^{\frac{1}{p}}, \frac{1}{4 \lambda\,2^{\frac{1}{p}}}\left(\frac{\lambda^p - \norm{U}^p}{\lambda^p - \norm{U}^p +1}\right)^{\frac{1}{p}}\right\} & \text{for $0<\norm{U}< \frac{1}{4 \lambda\,2^{\frac{1}{p}}}$}.
		\end{cases}
		$$
	\end{thm}

	In particular, we obtain the following lower bound of the $\lambda$-powered Bohr radius for operator-valued pluriharmonic functions whenever $U$ is the identity operator on $X=\mathcal{B}(\mathcal{H})$.
	\begin{cor}
		Let $X=\mathcal{B}(\mathcal{H})$ and $1<\lambda$. Then for all $n\in \mathbb{N}$ and $1\leq p < \infty$, we have
		\begin{equation*}
			R_{\lambda}(B_{Z}, p) \geq \frac{1}{2^{2+\frac{1}{p}}} \frac{\left(\lambda^p -1\right)^{\frac{1}{p}}}{\lambda^2}\, \frac{1}{\sup_{z \in B_{Z}}\norm{z}_{p}}.
		\end{equation*}
	\end{cor}
	\begin{rem}
		One might ask whether Theorem \ref{thm-1.1} remains valid when $\norm{U}=\lambda$. The following example shows that in this case the constant $R_{\lambda}(B_{Z}, p,U)$ may not be strictly positive. Thus the condition $\norm{U}<\lambda$ is indeed necessary.
	\end{rem}
	\begin{example} \label{example-1.1}
		For a bounded simply connected complete Reinhardt domain $\Omega \subseteq \mathbb{C}^n$ and $k \in \mathbb{N}$, define 
		$F_{k}: \Omega \rightarrow \mathcal{B}(\mathcal{H})$ by $$F_{k}(z)=(i\, \cos \frac{1}{k}) I_{\mathcal{H}} + (\frac{1}{2} \sin \frac{1}{k})  I_{\mathcal{H}}z_{1}+(\frac{1}{2} \sin \frac{1}{k})  I_{\mathcal{H}} \overline{z_{1}}$$
		 for $z=(z_{1}, \ldots,z_{n}) \in \Omega$, where $I_{\mathcal{H}}$ is the identity on $\mathcal{H}$. Let $U=\lambda \, I$, with $I$ the identity on $\mathcal{B}(\mathcal{H})$. Then $F_{k}(0)=(i\, \cos \frac{1}{k})  I_{\mathcal{H}}$. Clearly, $F_{k} \in \mathcal{PH}(\Omega,\mathcal{B}(\mathcal{H}))$ and $\norm{F_{k}}_{\Omega,\mathcal{B}(\mathcal{H})} \leq 1$. Assume, for contradiction, that there exists $r_{0}>0$ such that \eqref{e-1.4-a} holds for all $f \in \mathcal{PH}(\Omega,\mathcal{B}(\mathcal{H}))$ for all $z \in r_{0}\, \Omega$. Applying this to $F_{k}$ yields, $$\lambda^p \, |\cos \frac{1}{k}|^p + \lambda^p \, |\sin \frac{1}{k}|^p\, |z_{1}|^p \leq \lambda^p$$
		  all $z \in r_{0}\, \Omega$ and for all $k \in \mathbb{N}$. However, this inequality fails for sufficiently large $k$, since $\cos (1/k) \rightarrow 1$ and $\sin(1/k) \rightarrow 0$ as $k \rightarrow\infty$,while $z$
		can be chosen in $r_{0}\, \Omega$ such that $z_{1} \neq 0$. This contradiction shows that the condition $\norm{U}<\lambda$ in Theorem \ref{thm-1.1} is necessary.
	\end{example}
	We now turn to the second part of Question \ref{qsn-1.3}, which is the central focus of this paper, namely the asymptotic estimate of the constant $R_{\lambda}(\Omega, p,U)$. The exact asymptotic behavior of the $\lambda$-Bohr radius $K_{\lambda}(\mathbb{D}^n, 1,U)$ for holomorphic functions is known from \cite{defant-2012}. One of our main contributions of this paper is to extend these results by determining the precise asymptotic behavior of the powered $\lambda$-Bohr radius for operator valued pluriharmonic functions on any complete Reinhardt domains. The investigation separates naturally into two cases, according to whether $X$ is finite and infinite dimensional. We first address the finite dimensional case, where the following theorem describes the exact asymptotic behavior of the constant $R_{\lambda}(\Omega, p,U)$ when $\Omega=B_{Z}$ and $U$ is the identity operator on $X=\mathcal{B}(\mathcal{H})$. Recall that a Schauder basis $\{w_{k}\}$ of a Banach space $W$ is said to be $1$-unconditional if $\chi(\{w_k\})=1$ (see more details in Section $3$). We usually denote the canonical basis vectors of $Z=(\mathbb{C}^n, ||.||)$ by $e_{k}$, $k=1, \ldots,n$. 
	\begin{thm} \label{thm-1.2}
		Let $Z=(\mathbb{C}^n, ||.||)$ be a Banach space such that $\chi(\{e_k\}^n_{k=1})=1$. Let $X=\mathcal{B}(\mathcal{H})$ be a finite dimensional complex Banach space and $\lambda>1$. Then 
		\comment{\begin{equation}
				R^{m}_{\lambda}(B_{Z}, p) \geq	\begin{cases}
					E_{1}(X)\, \lambda^{\frac{1}{m}} \, \max \left\{\frac{1}{\sqrt{n}\, \norm{Id: \ell^n_{2} \rightarrow Z}}, \, \frac{1}{\left(\frac{m^m}{m!}\right)^{1/m}\, \norm{Id: Z \rightarrow \ell^n_{1}}}\right\} & \text{for $p=1$}, \\[3mm]
					2^{-\frac{p+5}{pm}}\,  \lambda^{\frac{1}{m}}\, \norm{I:\ell^n _2 \rightarrow Z}^{-\frac{2}{p}} \,\norm{I-\real (a_{0})}^{-\frac{2}{pm}} & \text{for $p \geq 2$}, \\[2mm]
					
					E_{2}(X)\, \lambda^{\frac{1}{m}} \,\max \left\{\frac{1}{(\sqrt{n})^{1-\theta}\, \norm{Id: \ell^n_{2} \rightarrow Z}}, \, \frac{\norm{I:\ell^n _2 \rightarrow Z}^{-\theta}}{\left(\frac{m^m}{m!}\right)^{(1-\theta)/m}\, \norm{Id: Z \rightarrow \ell^n_{1}}^{1-\theta}}\right\} & \text{for $1<p< 2$}, 
				\end{cases}
			\end{equation}
			where $\theta=(2(p-1))/p$. Moreover,}
		\begin{equation*}
			R_{\lambda}(B_{Z}, p,X) \geq	\begin{cases}
				E_{1}(X)\,  \max \left\{\frac{1}{\sqrt{n}\, \norm{Id: \ell^n_{2} \rightarrow Z}}, \, \frac{1}{e\, \norm{Id: Z \rightarrow \ell^n_{1}}}\right\}\, \left(\frac{\lambda^p -1}{2\lambda^p - 1}\right)^{\frac{1}{p}} & \text{for $p=1$}, \\[2mm]
				E_{3}\norm{Id:\ell^n _2 \rightarrow Z}^{-\frac{2}{p}}\, \left(\frac{\lambda^p -1}{2\lambda^p - 1}\right)^{\frac{1}{p}}\, & \text{for $p \geq 2$}, \\[2mm]
				
				E_{2}(X)\,  \max \left\{\frac{1}{(\sqrt{n})^{1-\theta} \norm{Id: \ell^n_{2} \rightarrow Z}},  \frac{\norm{Id:\ell^n _2 \rightarrow Z}^{-\theta}}{e^{(1-\theta)} \norm{Id: Z \rightarrow \ell^n_{1}}^{1-\theta}}\right\}\left(\frac{\lambda^p -1}{2\lambda^p - 1}\right)^{\frac{1}{p}} \hspace{-2mm} & \text{for $1<p< 2$}, 
			\end{cases}
		\end{equation*}
		where $\theta = (2(p-1))/p$. Here, $E_{1}(X)$, $E_{2}(X)$, and $E_{3}$ are positive constants: $E_{1}(X)$ depends only on $X$, $E_{2}(X)$ depends on both $X$ and $p$, while $E_{3}$ is independent of $X$ and depends only on $p$. Furthermore,
		for any $1\leq q \leq \infty$,
		\begin{equation*}
			K_{\lambda}(B_{Z}, p,X) \leq d\,\norm{Id: Z \rightarrow \ell^n_{q}}\, \norm{Id: \ell^n_{q} \rightarrow Z}\, \lambda^{\frac{2}{\log \, n}}\,n^{1-\frac{1}{p}}\, \left(\frac{\log \, n}{n}\right)^{1-\left(1/\min\{q, 2\}\right)}
		\end{equation*}
		for some constant $d>0$.
	\end{thm}
	In particular, for the Minkowski spaces $Z=\ell^n _q$, we obtain the following lower bound. The upper for this obtained in the proof of Theorem \ref{thm-1.2}.
	\begin{cor} \label{cor-1.4}
		Let $Z=(\mathbb{C}^n, ||.||)$ be a Banach space such that $\chi(\{e_k\}^n_{k=1})=1$.	Let $X=\mathcal{B}(\mathcal{H})$ be finite dimensional and $\lambda>1$. Then
		\begin{equation*}
			R_{\lambda}(B_{\ell^n _q}, p,X) \geq	\begin{cases}
				E'_{3}(X)\,\left(\frac{\lambda^p -1}{2\lambda^p - 1}\right)^{\frac{1}{p}} \, \left(\frac{\log n}{n}\right)^{ \left(1- \frac{1}{\min\{q,2\}}\right)} & \text{for $p=1$}, \\[2mm]
				E'_{2} \left(\frac{\lambda^p -1}{2\lambda^p - 1}\right)^{\frac{1}{p}}\, n^{-\frac{1}{p}} & \text{for $p \geq q$}, \\[2mm]
				
				E_{4}(X)\,  \left(\frac{\lambda^p -1}{2\lambda^p - 1}\right)^{\frac{1}{p}} \,n^{-\frac{p-1}{p(q-1)}}\, \left( \frac{\log n}{n}\right)^{\left(1- \frac{1}{\min\{q,2\}}\right)\, \frac{q-p}{p(q-1)}} & \text{for $1<p< q$}. 
			\end{cases}
		\end{equation*}
		Here, $E'_{3}(X)$, $E'_{2}$, and $E_{4}(X)$ are positive constants: $E'_{3}(X)$ depends only on $X$, $E_{4}(X)$ depends on both $X$ and $p$, while $E'_{2}$ is independent of $X$ and depends only on $p$
	\end{cor}
	\noindent In light of Theorem \ref{thm-1.2}, we deduce the following asymptotic estimates for Banach sequence spaces, which have wide applications, as discussed in Section $2$.
	\begin{thm} \label{thm-1.3-a}
		Let $Z$ be a Banach sequence space, $X=\mathcal{B}(\mathcal{H})$ be finite dimensional, and $\lambda>1$. For $n\in \mathbb{N}$, let $Z_{n}$ be the linear span of $\{e_{k}, \, k=1,\ldots,n\}$. Then
		\begin{enumerate}
			\item if $Z$ is a subset of $\ell_2$ we have \begin{equation*}
				R_{\lambda}(B_{Z_{n}}, p,X) \geq	\begin{cases}
					E_{1}(X)\,  \, \max \left\{\frac{1}{\sqrt{n}\, \norm{Id: \ell^n_{2} \rightarrow Z_{n}}}, \, \frac{1}{e\, \norm{\sum_{k=1}^{n}e^{*}_{k}}_{Z^{*}}}\right\}\, \left(\frac{\lambda^p -1}{2\lambda^p - 1}\right)^{\frac{1}{p}} & \text{for $p=1$}, \\[2mm]
					E_{3}\norm{Id:\ell^n _2 \rightarrow Z_{n}}^{-\frac{2}{p}}\, \left(\frac{\lambda^p -1}{2\lambda^p - 1}\right)^{\frac{1}{p}}\, & \text{for $p \geq 2$}, \\[2mm]
					
					E_{2}(X)\,  \max \left\{\frac{1}{(\sqrt{n})^{1-\theta} \norm{Id: \ell^n_{2} \rightarrow Z_{n}}}, \frac{\norm{Id:\ell^n _2 \rightarrow Z_{n}}^{-\theta}}{e^{(1-\theta)} \norm{\sum_{k=1}^{n}e^{*}_{k}}_{Z^{*}}^{1-\theta}}\right\}\left(\frac{\lambda^p -1}{2\lambda^p - 1}\right)^{\frac{1}{p}} \hspace{-3mm} & \text{for $1<p< 2$}. 
				\end{cases}
			\end{equation*}
			and 
			$R_{\lambda}(B_{Z_{n}}, p,X) \leq d\,  \lambda^{2/\log \, n}\,\norm{Id: \ell^n_{2} \rightarrow Z_{n}}\,n^{1-1/p}\, \sqrt{\log \, n/n};$
			\item if $Z$ is symmetric and $2$-convex, we have
			\begin{equation*}
				R_{\lambda}(B_{Z_{n}}, p,X) \geq	\begin{cases}
					E_{1}(X)\,  \, \max \left\{\frac{1}{\sqrt{n}}, \, \frac{1}{e\, \norm{\sum_{k=1}^{n}e^{*}_{k}}_{Z^{*}}}\right\}\, \left(\frac{\lambda^p -1}{2\lambda^p - 1}\right)^{\frac{1}{p}} & \text{for $p=1$}, \\[2mm]
					E_{3}\, \left(\frac{\lambda^p -1}{2\lambda^p - 1}\right)^{\frac{1}{p}}\, & \text{for $p \geq 2$}, \\[2mm]
					
					E_{2}(X)\,  \max \left\{\frac{1}{(\sqrt{n})^{1-\theta}}, \, \frac{1}{e^{(1-\theta)}\, \norm{\sum_{k=1}^{n}e^{*}_{k}}_{Z^{*}}^{1-\theta}}\right\}\left(\frac{\lambda^p -1}{2\lambda^p - 1}\right)^{\frac{1}{p}}  & \text{for $1<p< 2$}. 
				\end{cases}
			\end{equation*}
			and 
			$R_{\lambda}(B_{Z}, p,X) \leq d\,  \lambda^{2/\log \, n}\, \norm{\sum_{k=1}^{n}e^{*}_{k}}_{Z^{*}}\,\,n^{-1/p}\,\, \sqrt{\log \, n}.
			$
		\end{enumerate}
		The constants $E_{1}, E_{2}, E_{3}$, and $d$ are same as in Theorem \ref{thm-1.2}.
	\end{thm}
	To obtain the exact asymptotic behavior of the constant $R_{\lambda}(B_{Z}, p,X)$ in the case where $X$ is infinite dimensional, we shall make use of Lemma \ref{lem-3.5}. This lemma reduces the problem to studying the constant $R_{\lambda}(\Omega, p,X)$ for a complete Reinhardt domain $\Omega$, where the task of estimating bounds becomes considerably simpler. The most natural examples of such domains in $\mathbb{C}^n$ are the unit balls $B_{\ell^n _q}$ of the Minkowski spaces $\ell^n _q$, with $1\leq q \leq \infty$. We therefore first analyze $R_{\lambda}(\Omega, p,X)$ when $\Omega=B_{\ell^n _q}$. For finite dimensional $X$, Corollary \ref{cor-1.4} shows that the asymptotics of $R_{\lambda}(B_{\ell^n _q}, p)$ involve a logarithmic term. In contrast, in the infinite dimensional case, the following theorem demonstrates that this logarithmic factor disappears, and the asymptotic decay of $R_{\lambda}(B_{\ell^n _q}, p)$ is instead determined by the geometry of the Banach space $X$, namely its optimal cotype $Cot(X)$(see Section $2$ for the definition). 
	\begin{thm} \label{thm-1.3}
		Let $X=\mathcal{B}(\mathcal{H})$ be an infinite dimensional complex Banach space  and $\lambda>1$. Then 
		\begin{equation*}
			R_{\lambda}(B_{\ell^n _q}, p,X) \geq \Psi_{1}(n,\lambda,p,q):=	\begin{cases}
				
				E_{5}\,\dfrac{\left(\lambda^p -1\right)^{\frac{1}{p}}}{\lambda} & \text{for $p>\min \{Cot(X),q\}$}, \\[2mm]
				E_{6}\, \dfrac{\left(\lambda^p -1\right)^{\frac{1}{p}}}{\lambda}\,\, n^{{\frac{1}{\min\{t,q\}}}-\frac{1}{p}} & \text{for $p\leq \min \{Cot(X),q\}$},
			\end{cases}
		\end{equation*}
		if $X$ has a finite cotype $t$. Here $E_{5}>0$ and $E_{6}>0$ are constant independent of $X$ and $n$. 
		On the other hand, we have
		\begin{equation*}
			R_{\lambda}(B_{\ell^n _q}, p,X) \leq \Psi_{2}(n,\lambda,p,q):=	\begin{cases}	
				2\, \lambda\, n^{\frac{1}{q}\, - \,\frac{1}{p}} & \text{for $q\leq Cot(X)$}, \\[2mm]
				2\,\lambda\, n^{\frac{1}{Cot(X)}\, - \,\frac{1}{p}} & \text{for $q>Cot(X)$}.
			\end{cases}
		\end{equation*}
	\end{thm}

	We are now in a stage to give the bounds of the Bohr radius $R_{\lambda}(B_{Z}, p,X)$ for the unit ball $B_{Z}$ of $Z=(\mathbb{C}^n,||.||)$ whenever $X$ is a infinite dimensional complex Banach space. In view of Theorem \ref{thm-1.3} and Lemma \ref{lem-3.5}, we obtain the following result. Before stating the result, we introduce the following standard notation. For bounded complete Reinhardt domains $\Omega_{1}, \Omega_{2}\subset \mathbb{C}^n$, let $S(\Omega_1, \Omega_{2}):= \inf \left\{s>0 : \Omega_{1} \subset s \, \Omega_{2}\right\}$. Recall that if $Z$ and $W$ are Banach sequence spaces then $S(B_{Z},B_{W})=\norm{Id:Z \rightarrow W}$.
	\begin{thm} \label{thm-1.4}
		Let $Z=(\mathbb{C}^n, ||.||)$ be a Banach space such that $\chi(\{e_k\}^n_{k=1})=1$. Let $\mathcal{B}(\mathcal{H})$ be infinite dimensional, and $\Psi_{1}(n,\lambda,p,q)$ and $\Psi_{2}(n,\lambda,p,q)$ be as in Theorem \ref{thm-1.3}. For $\lambda>1$, 
		$$\frac{\Psi_{1}(n,\lambda,p,q)}{S(B_{Z},B_{\ell^n _q})S(B_{\ell^n _q},B_{Z})} \, \leq R_{\lambda}(B_{Z}, p,X) \leq S(B_{Z},B_{\ell^n _q})S(B_{\ell^n _q},B_{Z}) \, \Psi_{2}(n,\lambda,p,q).
		$$
	\end{thm}
	\begin{thm} \label{thm-1.6}
		Let $Z$ be a Banach sequence space and $X=\mathcal{B}(\mathcal{H})$ be infinite dimensional. For $n\in \mathbb{N}$, let $Z_{n}$ be the linear span of $\{e_{k}, \, k=1,\ldots,n\}$. Then for $\lambda>1$
		\begin{enumerate}
			\item if $Z$ is subset of $\ell_2$, we have 
			$$\frac{\Psi_{1}(n,\lambda,p,q)}{S(B_{\ell^n _q},B_{Z_{n}})} \, \leq R_{\lambda}(B_{Z_{n}}, p,X) \leq S(B_{\ell^n _2},B_{Z_{n}}) \, \Psi_{2}(n,\lambda,p,q);$$
			\item if $Z$ is symmetric and $2$-convex, we have $$\frac{\sqrt{n}\Psi_{1}(n,\lambda,p,q)}{\norm{\sum_{k=1}^{n}e^{*}_{k}}_{Z^{*}}} \, \leq R_{\lambda}(B_{Z_{n}}, p,X) \leq \frac{\norm{\sum_{k=1}^{n}e^{*}_{k}}_{Z^{*}}}{\sqrt{n}} \, \Psi_{2}(n,\lambda,p,q).
			$$	
		\end{enumerate}
	\end{thm}
\begin{rem}
We have thus far determined the asymptotic behavior of $R_{\lambda}(\Omega, p,X)$ in the case $\Omega=B_{Z}$, where $B_{Z}$ is the open unit ball of the Banach space $Z=(\mathbb{C}^n, ||.||)$ whose canonical basis is normalized $1$-unconditional. This framework also covers the bounded convex simply connected complete Reinhardt domains (see the correspondence preceding Theorem \ref{thm-1.1}). In light of Lemma \ref{lem-3.5}, Theorem \ref{thm-1.2}, and Theorem \ref{thm-1.4}, these asymptotic results extend to every bounded simply connected complete Reinhardt domain, without requiring convexity.
\end{rem}
	We now turn our attention to answering Question \ref{qsn-1.1}. To begin with, it is useful to point out some important differences between the study of Bohr’s theorem for vector-valued holomorphic functions and that for pluriharmonic functions. First, in the pluriharmonic setting we restrict ourselves to simply connected complete Reinhardt domains in $\mathbb{C}^n$. This restriction is necessary because the series expansion \eqref{e-1.3-a} is valid only on such domains. Second, in this framework we focus exclusively on $\mathcal{B}(\mathcal{H})$-valued pluriharmonic functions, and correspondingly restrict the operator $U$ to $\mathcal{B}(\mathcal{H})$. The reason is that the notion of complex conjugation is well defined in $\mathcal{B}(\mathcal{H})$, which is essential for our analysis. On the other hand, in the study of Bohr’s theorem for vector-valued holomorphic functions, the situation is somewhat simpler. Here, it suffices to consider complete Reinhardt domains (without the additional requirement of simply connectedness) because the domain of convergence of multivariable power series is complete Reinhardt domain and within its domain of convergence it defines a holomorphic function. Moreover, since complex conjugation does not play any role in this context, we are not constrained to $\mathcal{B}(\mathcal{H})$-valued functions. Instead, we can work with holomorphic functions taking values in any complex Banach space. Therefore, in addressing Question \ref{qsn-1.1}, we adopt the following framework: we consider bounded holomorphic functions with values in an arbitrary complex Banach space $X$, defined on complete Reinhardt domains in $\mathbb{C}^n$. In the following theorem, we first establish that the constant $K_{\lambda}(\Omega, p, U)$ is nonzero whenever $\Omega = B_{Z}$, which serves as the holomorphic analogue of Theorem \ref{thm-1.1}. By applying the same arguments as in the pluriharmonic case (see the discussion preceding Theorem \ref{thm-1.1}), it then follows that $K_{\lambda}(\Omega, p, U)$ is nonzero for every bounded complete Reinhardt domain in $\mathbb{C}^n$.
	\begin{thm} \label{thm-1.1-a}
		Let $X$ and $Y$ be any complex Banach spaces and $U: X \rightarrow Y$ be a non-null bounded linear operator such that $\norm{U}< \lambda$. Then, for $\lambda>1$ and $n \in \mathbb{N}$, we have
		$$K_{\lambda}(B_{Z}, p,U) \geq  \frac{C}{\sup_{z \in B_{Z}}\norm{z}_{p}},$$
		 where 
		$	C=\begin{cases}
			\max \left\{\left(\frac{\lambda^p - \norm{U}^p}{2\lambda^p - \norm{U}}\right)^{\frac{1}{p}}\, , \left(\frac{\lambda^p - \norm{U}^p}{\lambda^p - \norm{U}^p +1}\right)^{\frac{1}{p}}\, \frac{1}{\norm{U}}\right\}\,  & \text{for $\norm{U}\geq 1$},\\[3mm]
			\max \left\{\left(\frac{\lambda^p - \norm{U}^p}{2\lambda^p - \norm{U}}\right)^{\frac{1}{p}}, \left(\frac{\lambda^p - \norm{U}^p}{\lambda^p - \norm{U}^p +1}\right)^{\frac{1}{p}}\right\} & \text{for $0<\norm{U}< 1$}.
		\end{cases}
		$
	\end{thm}
	\noindent For $X=\mathcal{B}(\mathcal{H})$, the bounds of $R_{\lambda}(B_{Z}, p,X)$ in Theorems \ref{thm-1.2}–\ref{thm-1.6} extend to $K_{\lambda}(B_{Z}, p,X)$ for any complex Banach space $X$ with different constants, ensuring that the asymptotic behaviors of $R_{\lambda}(B_{Z}, p,X)$ and $K_{\lambda}(B_{Z}, p,X)$ coincide for $X=\mathcal{B}(\mathcal{H})$. The proofs are analogous and omitted.
	\section{Application}
	This section presents several applications of our results to various classes of sequence spaces. We begin with the family of mixed Minkowski spaces. Recall the definition of mixed Minkowski spaces: 
	$$\ell^m_s(\ell^n_t):=\{(z_{k})^m_{k=1}:z_{1}, \ldots,z_{m} \in \mathbb{C}^n\}$$ 
	with $\norm{(z_{k})^m_{k=1}}_{s,t}:=(\sum_{k=1}^{m}\norm{z_{k}}^s_t)^{1/s}$. When $s=t$, the spaces $\ell^m_s(\ell^n_t)$ reduces exactly to the classical Minkowski spaces $\ell^{mn}_s$. Earlier studies such as those in \cite{hamada-JFA-2022,defant-2012} focused primarily on classical Minkowski spaces. The framework developed in this article incorporates that setting as a particular case while at the same time treating a substantially broader family of sequence spaces. This unified approach enables us to study both holomorphic and pluriharmonic functions in general sequence spaces.
	\begin{example}[{\bf Mixed spaces}]
		Let $m,n \geq 2$. Then
		\begin{enumerate}
			\item Let $1 \leq s ,t \leq 2$. If $X=\mathcal{B}(\mathcal{H})$ is finite dimensional, then
			\begin{equation*}
				R_{\lambda}(B_{\ell^m_s(\ell^n_t)}, p,X) \geq	\begin{cases}
					E_{1}(X)\,  \, \max \left\{\frac{1}{n^{1/t}\,m^{1/s}}, \, \frac{1}{e\,n^{1-\frac{1}{t}}m^{1-\frac{1}{s}} }\right\}\, \left(\frac{\lambda^p -1}{2\lambda^p - 1}\right)^{\frac{1}{p}} & \text{for $p=1$}, \\[2mm]
					E_{3}\, n^{\frac{1}{p}-\frac{2}{tp}}\,m^{\frac{1}{p}-\frac{2}{sp}} \left(\frac{\lambda^p -1}{2\lambda^p - 1}\right)^{\frac{1}{p}}\, & \text{for $p \geq 2$}, \\[2mm]
					
					E_{2}(X)  \max \left\{\frac{1}{ n^{\frac{1}{t}+\frac{1}{p}-1}m^{\frac{1}{s}+\frac{1}{p}-1}},  \frac{m^{2+\frac{4}{ps}-\frac{3}{p}-\frac{3}{s}}}{e^{(1-\theta)} n^{\frac{3}{p}+\frac{3}{t} -2-\frac{4}{pt}}}\right\}\left(\frac{\lambda^p -1}{2\lambda^p - 1}\right)^{\frac{1}{p}} \hspace{-1mm} & \text{for $1<p< 2$}, 
				\end{cases}
			\end{equation*}
			and 
			$K_{\lambda}(B_{\ell^m_s(\ell^n_t)}, p,X) \leq d\, n^{\frac{1}{t}-\frac{1}{p}}m^{\frac{1}{s}-\frac{1}{p}}\, \lambda^{\frac{2}{\log \, n}}\, \sqrt{\log mn}.
			$
			If $X=\mathcal{B}(\mathcal{H})$ is infinite dimensional then
			$$\frac{\Psi_{1}(n,\lambda,p,2)}{m^{1/s -1/2}n^{1/t-1/2}} \, \leq R_{\lambda}(B_{Z}, p,X) \leq m^{1/s -1/2}n^{1/t-1/2} \, \Psi_{2}(n,\lambda,p,2).
			$$ 
			\item Let $2 \leq s,t \leq \infty$. If $X=\mathcal{B}(\mathcal{H})$ is finite dimensional then
			\begin{equation*}
				R_{\lambda}(B_{\ell^m_s(\ell^n_t)}, p,X) \geq	\begin{cases}
					E_{1}(X)\,  \, \max \left\{\frac{1}{\sqrt{mn}}, \, \frac{1}{e\, n^{1-\frac{1}{t}}m^{1-\frac{1}{s}}}\right\}\, \left(\frac{\lambda^p -1}{2\lambda^p - 1}\right)^{\frac{1}{p}} & \text{for $p=1$}, \\[2mm]
					E_{3}\, \left(\frac{\lambda^p -1}{2\lambda^p - 1}\right)^{\frac{1}{p}}\, & \text{for $p \geq 2$}, \\[2mm]
					
					E_{2}(X)  \max \left\{\frac{1}{(\sqrt{mn})^{1-\theta}}, \frac{1}{e^{(1-\theta)} n^{1-\theta -\frac{1-\theta}{t}}m^{1-\theta -\frac{1-\theta}{s}}}\right\}\left(\frac{\lambda^p -1}{2\lambda^p - 1}\right)^{\frac{1}{p}}  \hspace{-2mm} & \text{for $1<p< 2$}, 
				\end{cases}
			\end{equation*}
			and
			$K_{\lambda}(B_{\ell^m_s(\ell^n_t)}, p,X) \leq d\,n^{1-\frac{1}{t}-\frac{1}{p}}m^{1-\frac{1}{s}-\frac{1}{p}}\, \lambda^{\frac{2}{\log \, n}}\, \sqrt{\log mn}.
			$
			If $X=\mathcal{B}(\mathcal{H})$ is infinite dimensional then
			$$\frac{\Psi_{1}(n,\lambda,p,2)}{m^{1/2 -1/s}\,n^{1/2-1/t}} \, \leq R_{\lambda}(B_{Z}, p,X) \leq m^{1/2 -1/s}\,n^{1/2-1/t} \, \Psi_{2}(n,\lambda,p,2).
			$$
		\end{enumerate}
		Here, the constants $E_{1}(X)$, $E_{2}(X)$, $E_{3}$, and $d$ are as in Theorem \ref{thm-1.2}.	
	\end{example}
	\begin{pf}
		It is straightforward to verify that, for $1\leq s,t,l,r \leq \infty$, 
		$\norm{Id:\ell^m_s(\ell^n_t) \rightarrow \ell^m_l(\ell^n_r)}= \norm{Id:\ell^m_s \rightarrow \ell^m_l}\, \norm{Id:\ell^n_t \rightarrow \ell^n_r}$ (see \cite[p. 188]{defant-2003}). Furthermore, the norm of the identity operator between finite-dimensional $\ell^n_p$ spaces is given by $\norm{Id:\ell^n_s \rightarrow \ell^n_t}=1$ for $s \leq t$ and $\norm{Id:\ell^n_s \rightarrow \ell^n_t}=n^{1/t - 1/s}$ for $s>t$. Therefore, the desired result follows from Theorems \ref{thm-1.2}, \ref{thm-1.4}, and \ref{thm-1.6}, together with the preceding observations. The proof is complete.
	\end{pf}
	Our results can also be applied to many concrete symmetric Banach sequence spaces such as Lorentz spaces $\ell_{s,t}$ as well as Orlicz spaces $\ell_{\psi}$. For their definition we refer \cite{Lindenstrauss-book}. 
	\begin{example}[{\bf Lorentz spaces}]
		Let $1 \leq s,t \leq \infty$.
		\begin{enumerate}
			\item Let $1 \leq s <2$ and $1\leq t \leq \infty$, or $s=2$ and $1 \leq t \leq 2$. If $X=\mathcal{B}(\mathcal{H})$ is finite dimensional 
			\begin{equation*}
				R_{\lambda}(B_{\ell^n_{s,t}}, p,X) \geq	\begin{cases}
					E_{1}(X)\,  \, \max \left\{\frac{1}{\sqrt{n}\, \norm{Id: \ell^n_{2} \rightarrow \ell^n_{s,t}}}, \, \frac{1}{e\, n^{1-\frac{1}{s}}}\right\}\, \left(\frac{\lambda^p -1}{2\lambda^p - 1}\right)^{\frac{1}{p}} & \text{for $p=1$}, \\[2mm]
					E_{3}\norm{Id:\ell^n _2 \rightarrow \ell^n_{s,t}}^{-\frac{2}{p}}\, \left(\frac{\lambda^p -1}{2\lambda^p - 1}\right)^{\frac{1}{p}}\, & \text{for $p \geq 2$}, \\[2mm]
					
					E_{2}(X)  \max \left\{\frac{1}{(\sqrt{n})^{1-\theta} \norm{Id: \ell^n_{2} \rightarrow \ell^n_{s,t}}},  \frac{\norm{Id:\ell^n _2 \rightarrow \ell^n_{s,t}}^{-\theta}}{e^{(1-\theta)} n^{(1-\theta)\left(1-\frac{1}{s}\right)}}\right\}\left(\frac{\lambda^p -1}{2\lambda^p - 1}\right)^{\frac{1}{p}} \hspace{-2mm} & \text{for $1<p< 2$}, 
				\end{cases}
			\end{equation*}
			and 
			$R_{\lambda}(B_{\ell^n_{s,t}}, p,X) \leq d\,  \lambda^{2/\log \, n}\,\norm{Id: \ell^n_{2} \rightarrow \ell^n_{s,t}}\, n^{1/2-1/p}\, \sqrt{\log \, n}.
			$
			If $X=\mathcal{B}(\mathcal{H})$ is infinite dimensional, then
			$$\frac{\Psi_{1}(n,\lambda,p,2)}{\norm{Id:\ell^n _2 \rightarrow \ell^n_{s,t}}} \, \leq R_{\lambda}(B_{Z}, p,X) \leq \norm{Id:\ell^n _2 \rightarrow \ell^n_{s,t}} \, \Psi_{2}(n,\lambda,p,2).
			$$
			\item Let $2<s \leq \infty$ and $1 \leq t \leq \infty$. If $X=\mathcal{B}(\mathcal{H})$ is finite dimensional
			\begin{equation*} 
				R_{\lambda}(B_{\ell^n_{s,t}}, p,X) \geq	\begin{cases}
					E_{1}(X)\,  \, \max \left\{\frac{1}{\sqrt{n}}, \, \frac{1}{e\, n^{\left(1-\frac{1}{s}\right)}}\right\}\, \left(\frac{\lambda^p -1}{2\lambda^p - 1}\right)^{\frac{1}{p}} & \text{for $p=1$}, \\[2mm]
					E_{3}\, \left(\frac{\lambda^p -1}{2\lambda^p - 1}\right)^{\frac{1}{p}}\, & \text{for $p \geq 2$}, \\[2mm]
					
					E_{2}(X)\,  \max \left\{\frac{1}{(\sqrt{n})^{1-\theta}}, \, \frac{1}{e^{(1-\theta)}\, n^{(1-\theta)\left(1-\frac{1}{s}\right)}}\right\}\left(\frac{\lambda^p -1}{2\lambda^p - 1}\right)^{\frac{1}{p}}  & \text{for $1<p< 2$}, 
				\end{cases}
			\end{equation*}
			and $R_{\lambda}(B_{\ell^n_{s,t}}, p,X) \leq d\,  \lambda^{2/\log \, n}\, n^{\left(1-\frac{1}{p}-\frac{1}{s}\right)}\, \sqrt{\log \, n}.
			$ If $X=\mathcal{B}(\mathcal{H})$ is infinite dimensional, then $$\frac{\Psi_{1}(n,\lambda,p,2)}{\norm{Id:\ell^n_{s,t} \rightarrow \ell^n _2}} \, \leq R_{\lambda}(B_{Z}, p,X) \leq \norm{Id:\ell^n_{s,t} \rightarrow \ell^n _2} \, \Psi_{2}(n,\lambda,p,2).
			$$
		\end{enumerate}
	\end{example}
	\begin{pf}
		We first note that $\norm{\sum_{k=1}^{n}e^{*}_{k}}_{\ell_{s,t}^{*}}=n^{1-1/s}$ and the Lorentz spaces $\ell_{s,t}$ admit the following lexicographical inclusion order: $\ell_{p,q} \subseteq \ell_{s,t}$ if and only if $p<s$, or $p=s$ and $q \leq t$. Hence, part $(1)$ follows directly from these observations and Theorems \ref{thm-1.3-a}$(1)$, \ref{thm-1.4}, and \ref{thm-1.6}. For part $(2)$, consider first $2<s \leq \infty$ and $2 \leq t \leq \infty$. In this case, $\ell_{s,t}$ is $2$-convex (see \cite[p.189]{defant-2003}), and hence the assertion follows from Theorems \ref{thm-1.3-a}$(2)$, \ref{thm-1.4}, and \ref{thm-1.6}. For the remaining case $2 < s \leq \infty$, $1 \leq t \leq 2$, we invoke Theorems \ref{thm-1.2}, \ref{thm-1.4}, and \ref{thm-1.6}, together with the estimates $\norm{Id: \ell^n_{2} \rightarrow \ell^n_{s,t}} \leq 1$ and $\norm{Id: \ell^n_{s,t} \rightarrow \ell^n_{2}} \leq n^{1/2 - 1/s}$ (see \cite[p.189]{defant-2003}). This completes the proof.
	\end{pf}
	\begin{example} [{\bf Orlicz spaces}]
		Let $\psi$ be an Orlicz function which satisfies the $\Delta_{2}$-condition. 
		\begin{enumerate}
			\item Let $a^2 \leq T \psi(a)$ for all $a$ and some $T>0$.  If $X=\mathcal{B}(\mathcal{H})$ is finite dimensional, then
			\begin{equation*}
				R_{\lambda}(B_{\ell^n_{\psi}}, p,X) \geq	\begin{cases}
					E_{1}(X)\,  \, \max \left\{\frac{1}{\sqrt{n}\, \norm{Id: \ell^n_{2} \rightarrow \ell^n_{\psi}}}, \, \frac{1}{e\,n\, \psi^{-1}(1/n)}\right\}\, \left(\frac{\lambda^p -1}{2\lambda^p - 1}\right)^{\frac{1}{p}} & \text{for $p=1$}, \\[2mm]
					E_{3}\norm{Id:\ell^n _2 \rightarrow \ell^n_{\psi}}^{-\frac{2}{p}}\, \left(\frac{\lambda^p -1}{2\lambda^p - 1}\right)^{\frac{1}{p}}\, & \text{for $p \geq 2$}, \\[2mm]
					
					E_{2}(X)  \max \left\{\frac{1}{(\sqrt{n})^{1-\theta} \norm{Id: \ell^n_{2} \rightarrow \ell^n_{\psi}}},  \frac{\norm{Id:\ell^n _2 \rightarrow \ell^n_{\psi}}^{-\theta}}{(en)^{(1-\theta)} (\psi^{-1}(1/n))^{1-\theta}}\right\}\left(\frac{\lambda^p -1}{2\lambda^p - 1}\right)^{\frac{1}{p}} \hspace{-4mm}  & \text{for $1<p< 2$}, 
				\end{cases}
			\end{equation*}
			and $R_{\lambda}(B_{\ell^n_{\psi}}, p,X) \leq d\,  \lambda^{2/\log \, n}\,\norm{Id: \ell^n_{2} \rightarrow \ell^n_{\psi}}\,n^{1/2-1/p}\, \sqrt{\log \, n}.
			$
			If $X=\mathcal{B}(\mathcal{H})$ is infinite dimensional, 
			$$\frac{\Psi_{1}(n,\lambda,p,2)}{\norm{Id:\ell^n _2 \rightarrow \ell^n_{\psi}}} \, \leq R_{\lambda}(B_{\ell^n_{\psi}}, p,X) \leq \norm{Id:\ell^n _2 \rightarrow \ell^n_{\psi}} \, \Psi_{2}(n,\lambda,p,2).
			$$
			
			\item Let $\psi(\beta \,a) \leq T \beta^2\,\psi(a)$ for $0 \leq \beta,a \leq 1$ and some $T>0$. If $X=\mathcal{B}(\mathcal{H})$ is finite dimensional,
			\begin{equation*}
				R_{\lambda}(B_{\ell^n_{\psi}}, p,X) \geq	\begin{cases}
					E_{1}(X)\,  \, \max \left\{\frac{1}{\sqrt{n}}, \, \frac{1}{e\,n\, \psi^{-1}(1/n)}\right\}\, \left(\frac{\lambda^p -1}{2\lambda^p - 1}\right)^{\frac{1}{p}} & \text{for $p=1$}, \\[2mm]
					E_{3}\, \left(\frac{\lambda^p -1}{2\lambda^p - 1}\right)^{\frac{1}{p}}\, & \text{for $p \geq 2$}, \\[2mm]
					
					E_{2}(X)\,  \max \left\{\frac{1}{(\sqrt{n})^{1-\theta}}, \, \frac{1}{(en)^{(1-\theta)}\, (\psi^{-1}(1/n))^{1-\theta}}\right\}\left(\frac{\lambda^p -1}{2\lambda^p - 1}\right)^{\frac{1}{p}}  & \text{for $1<p< 2$}, 
				\end{cases}
			\end{equation*}
			and 
			$R_{\lambda}(B_{\ell^n_{\psi}}, p,X) \leq d\,  \lambda^{2/\log \, n}\, n^{1-1/p}\,\psi^{-1}(1/n)\,\, \sqrt{\log \, n}$.
			If $X=\mathcal{B}(\mathcal{H})$ is infinite dimensional, 
			$$\frac{\Psi_{1}(n,\lambda,p,2)}{\sqrt{n}\, \psi^{-1}(1/n)} \, \leq R_{\lambda}(B_{\ell^n_{\psi}}, p,X) \leq \sqrt{n}\, \psi^{-1}(1/n) \, \Psi_{2}(n,\lambda,p,2).
			$$
		\end{enumerate}
	\end{example}
	\begin{pf}
		To prove this we first observe the fact $\norm{\sum_{k=1}^{n}e^{*}_{k}}_{\ell_{\psi}^{*}}\, \norm{\sum_{k=1}^{n}e_{k}}_{\ell_{\psi}}=n$ (see \cite[p.192]{defant-2003}). The condition on $\psi$ in $(1)$ ensures $\ell_{\psi} \subseteq \ell_2$, while that in $(2)$ implies that $\ell_{\psi}$ is $2$-convex (see \cite[p.189]{defant-2003}). The result then follows from Theorems \ref{thm-1.3-a} and \ref{thm-1.6}, together with the preceding observations. This completes the proof.
	\end{pf}	
	\section{Preliminaries}
	We use standard notation and notions from Banach space theory. All Banach spaces $W$ are assumed to be complex, their topological duals are denoted by $W^{*}$ and their open unit balls by $B_{W}$. A Banach space $W$ is said to have cotype $t\in [0,\infty]$ if there exists a constant $C>0$ such that for arbitrarily chosen vectors $w_{1},\ldots,w_{n}\in W$ we have 
	$$\left(\sum_{k=1}^{n}\norm{w_{k}}^t\right)^{\frac{1}{t}} \leq C \left(\int_{0}^{1}\norm{\sum_{k=1}^{n}r_{n}(s)w_{k}}^2\, ds\right)^{\frac{1}{2}},$$
	where $r_{n}$ it the $n$-th Rademacher function on $[0,1]$ and we write 
	$Cot(W):=\inf\{2 \leq t \leq \infty: W\, \mbox{has cotype}\, t\}.
	$
	It is worth mentioning that every Banach space $W$ has cotype $\infty$, and whenever $Cot(W)=\infty$ we denote $1/Cot(W)=0$. A Schauder basis $\{w_{k}\}$ of a Banach space $W$ is said to be unconditional if there is a constant $c\geq 1$ such that $\norm{\sum_{j=1}^{k}\epsilon_{j}\mu_{j}w_{j}} \leq c \norm{\sum_{j=1}^{k}\mu_{j}w_{j}}
	$
	 for all $k \in \mathbb{N}$, all $\mu _{j} \in \mathbb{C}$, and all $\epsilon_{j} \in \mathbb{C}$ with $|\epsilon_{j}| \leq 1$, and in this case the best constant is denoted by $\chi(\{w_{k}\})$, called unconditional basis constant of $\{w_{k}\}$. We say that $\{w_{k}\}$ is $1$-unconditional basis whenever $\chi(\{w_{k}\})=1$. Moreover, the unconditional basis constant of $W$ is defined to be $\chi(W):=\inf \chi(\{w_{k}\})\in [1,\infty]$, the infimum taken over all unconditional basis $\{w_{k}\}$ of $W$. A Banach space $Z$ for which $\ell_1 \subset Z \subset c_{0}$ (with normal inclusions) is said to be a Banach sequence space if the canonical sequence $\{e_{k}\}$ forms a $1$-unconditional basis. Recall that a Banach lattice $Z$ is said to be $2$-convex if there exists a constant $C>0$ such that 
	 $$\norm{\left(\sum_{k=1}^{n}|x_{k}|^2\right)^{\frac{1}{2}}} \leq C 	\left(\sum_{k=1}^{n}\norm{x_{k}}^2\right)^{\frac{1}{2}}
	 $$
	  for all $x_{1}, \ldots, x_{k} \in Z$. Let $Z=(\mathbb{C}^n, ||.||)$ be a Banach space and $Y$ be any Banach space, and $m \in \mathbb{N}$. We denote by $\mathcal{P}(^m Z,Y)$ the space of all $m$-homogeneous polynomials $Q:Z \rightarrow Y$ of the form $Q(z)=\sum_{|\alpha|=m}c_{\alpha}z^{\alpha}$, together with the norm $\norm{Q}_{\mathcal{P}(^m Z,Y)}:=\sup_{z \in B_{Z}}|Q(z)|$. By $\chi_{M}(\mathcal{P}(^m Z,Y))$ we denote the unconditional basis constant of the monomials $z^{\alpha}$, $\alpha \in (\mathbb{N} \cup \{0\})^n$.
	Set $X:=\mathcal{B}(H)$, the Banach space of all bounded linear operators on a Hilbert space $H$. Denote by $\mathcal{PH}(^m Z,X)$ the space of all $m$-homogeneous pluriharmonic polynomials $P:Z \rightarrow X$ of the form 
	\begin{equation} \label{e-1.3-aa}
		P(z)= \sum_{|\alpha|=m} a_{\alpha}\, z^{\alpha} +  \sum_{|\alpha|=m} b^{*}_{\alpha}\, \bar{z}^{\alpha}
	\end{equation} 
	and set $\norm{P}_{\mathcal{PH}(^m Z,X)}:=\sup_{z \in B_{Z}}|||P(z)||$. 
	\par A bounded linear operator $U:W\rightarrow Y$ between two Banach spaces $W$ and $Y$ is called $(r,t)$-summing, $r,t \in [1,\infty)$, if there is a constant $C>0$ such that for each choice of finitely many $w_{1}, \ldots,w_{n} \in W$ we have that 
	$$\left(\sum_{k=1}^{n} \norm{U(w_{k})}^r\right)^{1/r} \leq C \, \sup_{\phi \in B_{X^{*}}} \left(\sum_{k=1}^{n} |\phi(w_{k})|^t\right)^{1/t}.
	$$
	The best such $C$ is usually denoted by $\pi_{r,t}(U)$. In the case $r=t$ we call $U$ is $r$- summing and write the best constant as $\pi_{r}(U)$.
	\section{Key Properties and Supporting Lemmas} 
	In the scalar case, a standard approach to deriving nontrivial estimates for Bohr radii of holomorphic functions is to first analyze $m$-homogeneous polynomials, as this reduction often isolates the essential difficulties while providing sharper insight into the general framework. The following definition introduces the $m$-homogeneous analogue of Definition \ref{def-1.1}, which will serve as the basis for our further developments. Let $\Omega\subset \mathbb{C}^n$ be a simply connected complete Reinhardt domain and $n\in \mathbb{N}$. Let $U:X\rightarrow Y$ be a bounded liner operator and $\norm{U} \leq \lambda$. For $1 \leq p < \infty$, the $\lambda$-powered Bohr radius of $U$ with respect to $\mathcal{PH}(^m \Omega,X)$, denoted by $R^m_{\lambda}(\Omega, p,U)$, is defined to be the supremum of all $r\geq 0$ such that for every pluriharmonic polynomials $P:\Omega \rightarrow X$ of the form \eqref{e-1.3-aa} we have 
	\begin{equation} \label{e-1.4-aa}
		\sup_{z \in r\Omega}\, \sum_{|\alpha|=m} (\norm{U(a_{\alpha})}^p_{Y} + \norm{U(b_{\alpha})}^p_{Y})|z^\alpha|^p \leq \lambda^p\,\norm{f}^p_{\Omega,X}. 
	\end{equation}
	
	\noindent Here $\norm{f}_{\Omega,X}:=\sup_{z \in \Omega}\,\norm{f(z)}_{X}$. Set $R^m(\Omega, p,U):=R^m_{1}(\Omega, p,U)$, $R^m_{\lambda}(\Omega, p,X):=R^m_{\lambda}(\Omega, p,U)$ whenever $U=I:X\rightarrow X$, $R^m(\Omega, p,X):=R^m_{1}(\Omega, p,X)$, $R^m_{\lambda}(\Omega, p):=R^m_{\lambda}(\Omega, p,\mathbb{C})$, and $R^m(\Omega,p):=R^m_1(\Omega,p)$. It is straightforward that $R^m_{\lambda}(\Omega, p,U)=\lambda^{1/m}\,R^m(\Omega, p,U)$. Similarly as above, we can define the constant $K^m_{\lambda}(\Omega, p,U)$ for the class $\mathcal{P}^n_m(W)$ {\it i.e.}, for any Banach space $W$ valued $m$-homogeneous holomorphic polynomial $Q(z)=\sum_{|\alpha|=m}c_{\alpha}\, z^{\alpha}$. Clearly, we have $K^m_{\lambda}(\Omega, p,U)=\lambda^{1/m}\,K^m(\Omega, p,U)$. In light of Bohnenblust-Hille inequality, it is known that the constant $K^m(\mathbb{D}^n,1)$ is hypercontractive due to Defant {\it et al.} \cite{defant-2011}.
	\vspace{2mm}
	
	Schwarz–Pick type lemmas and distortion estimates for holomorphic and harmonic functions constitute an important direction of research in geometric function theory. Coefficient-type Schwarz–Pick estimates for complex-valued harmonic functions on the unit disc were obtained in \cite{Abu-2010, chen-2011}, and subsequently extended to complex-valued pluriharmonic functions on $B_{\ell^n_q}$ in \cite[Theorem 2.7]{hamada-JFA-2022}. 
	In analogy with \cite[Theorem 2.7]{hamada-JFA-2022}, we now extend these estimates to the operator valued setting by establishing a coefficient-type Schwarz–Pick lemma for pluriharmonic functions defined on general complete Reinhardt domains. This lemma constitutes a fundamental ingredient in the proofs of several main results of this article.
	\begin{lem} \label{lem-3.1} 
		Let $f \in \mathcal{PH}(Z,X)$ be of the form \eqref{e-1.3-a}. Then for all $|\alpha|=m\geq 1$, we have $\norm{\sum_{|\alpha|=m}(a_{\alpha} \pm b_{\alpha})\, z^{\alpha}}_{B_{Z},X} \leq 4 \,\norm{\norm{f}_{B_{Z},X}\,I-\real(a_{0})}$. Moreover, if $B_{Z}=B_{\ell^n_{q}}$ then $\norm{a_{\alpha}+b_{\alpha}}	\leq \frac{4}{\pi} \, \rho_{\alpha}\, \norm{\sum_{|\alpha|=m}(a_{\alpha}+b_{\alpha})\, z^{\alpha}}_{B_{Z},X}$,  $\norm{a_{\alpha}-b_{\alpha}}	\leq \frac{4}{\pi}\,\rho_{\alpha}\, \norm{\sum_{|\alpha|=m}(a_{\alpha}-b_{\alpha})\, z^{\alpha}}_{B_{Z},X}$, where $\rho_{\alpha}:=\left(\frac{|\alpha|^{|\alpha|}}{\alpha^{\alpha}}\right)^{1/q}$.
		\comment{\begin{enumerate}
				\item $\norm{\sum_{|\alpha|=m}(a_{\alpha}+b_{\alpha})\, z^{\alpha}}_{B_{Z},X} \leq 4 \,\gamma_{0} \norm{f}_{B_{Z},X}$,
				$\norm{\sum_{|\alpha|=m}(a_{\alpha}-b_{\alpha})\, z^{\alpha}}_{B_{Z},X} \leq 4 \, \gamma_{0}\, \norm{f}_{B_{Z},X}$;
				\item $\norm{a_{\alpha}+b_{\alpha}}	\leq \frac{4}{\pi} \, \rho_{\alpha}\, \norm{\sum_{|\alpha|=m}(a_{\alpha}+b_{\alpha})\, z^{\alpha}}_{B_{Z},X}$,  $\norm{a_{\alpha}-b_{\alpha}}	\leq \frac{4}{\pi}\,\rho_{\alpha}\, \norm{\sum_{|\alpha|=m}(a_{\alpha}-b_{\alpha})\, z^{\alpha}}_{B_{Z},X}$,
			\end{enumerate}
			where $\gamma_{0}:=\norm{I-\real(a_{0})}$ and $\rho_{\alpha}:=\left(\frac{|\alpha|^{|\alpha|}}{\alpha^{\alpha}}\right)^{\frac{1}{q}}$.}
	\end{lem}
	\begin{pf}
		For any fixed $z\in B_{Z}$, we consider the function $g: \mathbb{D} \rightarrow \mathcal{B}(H)$ by $g(\xi)=f(\xi z)$ for $\xi \in \mathbb{D}$. Then $g$ is a $\mathcal{B}(H)$-valued harmonic function on $\mathbb{D}$ with $\norm{g}_{\mathbb{D},X} \leq \norm{f}_{B_{Z},X}$. Observe that $g$ has the following expansion 
		$$g(\xi)= \sum_{m=0}^{\infty} \left(\sum_{|\alpha|=m} a_{\alpha}z^{\alpha}\right)\, \xi^m + \sum_{m=1}^{\infty} \left(\sum_{|\alpha|=m} b_{\alpha}z^{\alpha}\right)^{*}\, \bar{\xi}^{m}$$
		for $\xi \in \mathbb{D}$, and hence from the \cite[(2.12)]{bhowmik-2021} we conclude $\norm{\sum_{|\alpha|=m}(a_{\alpha}+b_{\alpha})\, z^{\alpha}}_{B_{Z},X} \leq 4 \norm{\norm{g}_{\mathbb{D},X}\,I-\real(a_{0})} \leq 4 \norm{\norm{f}_{B_{Z},X}\,I-\real(a_{0})}$.
		As $h=-if \in \mathcal{PH}(Z,X)$ and $\norm{h}_{B_{Z},X}=\norm{f}_{B_{Z},X}$, the desired estimate for $\norm{\sum_{|\alpha|=m}(a_{\alpha}-b_{\alpha})\, z^{\alpha}}_{B_{Z},X}$ follows by applying the preceding inequality to $h$. Now assume that $B_{Z}=B_{\ell^n_{q}}$. We consider the operator-valued holomorphic function $\tilde{h}(z)=\sum_{|\alpha|=m} (a_{\alpha} + b_{\alpha})z^{\alpha}$ on $B_{\ell^n _q}$. For any $\phi \in X^{*}$ with $\norm{\phi} \leq 1$, $\tilde{g}=\phi \circ \tilde{h}$ is a complex-valued holomorphic function on $B_{\ell^n _q}$ with $\norm{\tilde{g}}_{B_{\ell^n_{q}},X} \leq \norm{\tilde{h}}_{B_{\ell^n_{q}},X}$. Moreover, $\tilde{g}(z)=\sum_{|\alpha|=m} \phi (a_{\alpha} +b_{\alpha})z^{\alpha}$ for $z \in B_{\ell^n _q}$. Observe that $\tilde{g}$ is also a complex-valued pluriharmonic function, and hence from \cite[Theorem 2.6]{hamada-JFA-2022}, we have for $|\alpha|=m\geq 1$,
		\begin{equation*} \label{e-3.5-a}
			\norm{a_{\alpha} +b_{\alpha}}=\sup _{\phi \in B_{X^{*}}} |\phi(a_{\alpha}+b_{\alpha})| \leq \frac{4}{\pi} \left(\frac{|\alpha|^{|\alpha|}}{\alpha^{\alpha}}\right)^{\frac{1}{q}} \norm{\tilde{g}}_{B_{\ell^n_{q}},X} \leq \frac{4}{\pi} \left(\frac{|\alpha|^{|\alpha|}}{\alpha^{\alpha}}\right)^{\frac{1}{q}} \norm{\tilde{h}}_{B_{\ell^n_{q}},X}.
		\end{equation*}
		Finally, the bound of $\norm{a_{\alpha} -b_{\alpha}}$ follows by considering the function $\tilde{P}=-i\tilde{h}$.
	\end{pf}

	The following result constitutes the basic link between the Bohr radius for pluriharmonic functions and the Bohr radius for homogeneous pluriharmonic functions. It plays a crucial role in the proofs of our main results. This lemma can be regarded as the pluriharmonic analogue of \cite[Lemma 3.2]{defant-2012}, which concerns the Minkowski spaces in the context of the powered Bohr inequality. The proof follows along similar lines as that of \cite[Lemma 3.2]{defant-2012}, and therefore, we omit the details.
	\begin{lem} \label{lem-3.3}
		Let $X=\mathcal{B}(\mathcal{H})$ and $Y$ any complex Banach space, and $U: X \rightarrow Y$ be a non-null bounded linear operator such that $\norm{U}< \lambda$. Then, for all $p\in[1,\infty), q \in[1,\infty]$, $\lambda>1$, and $n \in \mathbb{N}$, we have
		\begin{enumerate}
			\item $\left(\frac{\lambda^p - \norm{U}^p}{2\lambda^p - \norm{U}^p}\right)^{\frac{1}{p}}\, \inf_{m\in \mathbb{N}}\, R^{m}_{\lambda}(B_{Z}, p,U) \leq R_{\lambda}(B_{Z}, p,U) \leq \inf_{m\in \mathbb{N}}\, R^{m}_{\lambda}(B_{Z}, p,U)$
			
			\item $\left(\frac{\lambda^p - \norm{U}^p}{\lambda^p - \norm{U}^p +1}\right)^{\frac{1}{p}}\, \inf_{m\in \mathbb{N}}\, R^{m}(B_{Z}, p,U) \leq R_{\lambda}(B_{Z}, p,U) \leq \lambda \, \inf_{m\in \mathbb{N}}\, R^{m}(B_{Z}, p,U)$.
		\end{enumerate}
	\end{lem}
	\begin{rem} \label{rem-4.1}
		Analogous results of Lemma \ref{lem-3.3} hold for $K^m_{\lambda}(B_{Z}, p,U)$ and $K_{\lambda}(B_{Z}, p,U)$.
	\end{rem}
	\noindent The following lemma compares the Bohr radius constants for two complete Reinhardt domains. Its proof is similar to that of \cite[Lemma 2.5]{defant-2004} and is therefore omitted. 
	\begin{lem} \label{lem-3.5}
		Let $\Omega_{1}$ and $\Omega_{2}$ be two bounded complete Reinhardt domains in $\mathbb{C}^n$. Then we have $$\frac{R_{\lambda}(\Omega_{2}, p,U)}{S(\Omega_{1},\Omega_{2})S(\Omega_{2},\Omega_{1})} \, \leq R_{\lambda}(\Omega_{1}, p,U) \leq S(\Omega_{1},\Omega_{2})S(\Omega_{2},\Omega_{1}) \, R_{\lambda}(\Omega_{2}, p,U).
		$$
	\end{lem}
	\noindent As an immediate application, we observe that the Bohr radius $R_{\lambda}(\mathbb{D}^n, p,U)$ serves as a universal lower bound for the Bohr radius $R_{\lambda}(\Omega, p,U)$ of every complete Reinhardt domain.
	\begin{lem} \label{lem-3.4}
		For any complete Reinhardt domain $\Omega$, we have 
		$
		R_{\lambda}(\Omega, p,U) \geq R_{\lambda}(\mathbb{D}^n, p,U).
		$
	\end{lem}
	We use the following remarkable result by Maurey and Pisier to prove the desired upper bound in the case of infintedimensional complex Banach space. 
	\begin{customthm}{A} \cite[Theorem 14.5]{diestel-abs-summing-1995} \label{thm-A}
		Given any infinite dimensional complex Banach space $X$, there exist $x_{1}, \ldots , x_{n} \in X$ for each $n \in \mathbb{N}$ such that $\norm{z}_{\infty}/2 \leq \norm{\sum_{j=1}^{n} x_{j}z_{j}} \leq \norm{z}_{Cot(X)}$  for every choice of $z=(z_{1}, \ldots,z_{n}) \in \mathbb{C}^n$. Clearly, setting $z=e_{j}$, gives $\norm{x_{j}}\geq 1/2$, where $e_{j}$ is the $j$-th canonical basis vector of $\mathbb{C}^n$.
	\end{customthm}
	\section{Proof of the main results}
	We now begin the proofs of the main results. Let $f(z)=\sum_{|\alpha|=m} c_{\alpha}\, z^{\alpha}$ be an $X$ valued function. For any $p \geq 1$, we write $M_{f,p}:=\sum_{|\alpha|=m}\norm{c_{\alpha}\, z^{\alpha}}^p$, which will be used throughout the proofs.
	\begin{proof}[{\bf Proof of Theorem \ref{thm-1.1}}]
		Let $P(z)=\sum_{|\alpha|=m} a_{\alpha}\, z^{\alpha} + \sum_{|\beta|=m} b^{*}_{\alpha}\, \bar{z}^{\alpha} \in \mathcal{PH}(^m Z,X)$. Consider the polynomial $H(z):=\sum_{|\alpha|=m} \left(\frac{a_{\alpha} + b_{\alpha}}{2}\right)z^{\alpha}$. Applying the definition of the $\lambda$-powered Bohr radius for holomorphic functions to the function $H$  and then from Lemma \ref{lem-3.1}, for $z \in B_{Z}$, we obtain 
		$$M_{H,p} \leq \frac{\lambda^p}{\left(K^m_{\lambda}(B_{Z}, p,U)\right)^{pm}} \, \norm{H}^p_{B_{Z},X}\leq \frac{2^p\lambda^p}{\left(K^m_{\lambda}(B_{Z}, p,U)\right)^{pm}} \norm{P}^p_{B_{Z},X}.
		$$
		Similarly, considering $G(z):=\sum_{|\alpha|=m} \left(\frac{a_{\alpha} - b_{\alpha}}{2}\right)z^{\alpha}$,
		we obtain 
		$$M_{G,p}\leq \frac{2^p\lambda^p}{\left(K^m_{\lambda}(B_{Z}, p,U)\right)^{pm}} \norm{P}^p_{B_{Z},X}$$
			%
		for $z \in B_{Z}$. Now, for $p\geq 1$, using these inequalities, for $z \in B_{Z}$ we have
		\begin{align} \label{e-1.5-a}
			\left(\sum_{|\alpha|=m} \norm{a_{\alpha}z^{\alpha}}^p\right)^{1/p} \leq \left(M_{H,p}\right)^{1/p} + \left(M_{G,p}\right)^{1/p}
			\leq \frac{4\lambda}{\left(K^m_{\lambda}(B_{Z}, p,U)\right)^{m}}  \norm{P}_{B_{Z},X}.
		\end{align}
		
		\noindent One can verify that the inequality \eqref{e-1.5-a} remains true when the left-hand side is replaced by $\left(\sum_{|\alpha|=m} \norm{b_{\alpha}z^{\alpha}}^p\right)^{1/p}$ for $z \in B_{Z}$.
		By making use of these inequalities, for $z \in B_{Z}$,
		we obtain
		\begin{align}
			&\sum_{|\alpha|=m} \left(\norm{U(a_{\alpha})}^p + \norm{U(b_{\alpha})}^p\right) |z^{\alpha}|^p \nonumber
			\leq \left(\frac{K^m_{\lambda}(B_{Z}, p,U)}{2^{\frac{1}{pm}}(4\lambda)^{\frac{1}{m}}}\right)^{-pm}\lambda^p\norm{P}^p_{B_{Z},X},
		\end{align}
		which immediately shows that  
		\begin{equation} \label{e-1.8}
			R^m_{\lambda}(B_{Z}, p,U) \geq \frac{K^m_{\lambda}(B_{Z}, p,U)}{2^{\frac{1}{pm}}(4\lambda)^{\frac{1}{m}}}.
		\end{equation}
		According to \eqref{e-1.8}, it is enough to focus on finding a lower bound for $K^m_{\lambda}(B_{Z}, p,U)$, as this directly leads to a lower bound for $R^m_{\lambda}(B_{Z}, p,U)$. Let $Q(z)=\sum_{|\alpha|=m} a_{\alpha}\, z^{\alpha}  \in \mathcal{P}(^m Z,X)$. For any $z \in B_{Z}$ we have $$\sum_{|\alpha|=m} \norm{U(a_{\alpha})z^{\alpha}}^p \leq  \lambda^p  \norm{Q}^p_{B_{Z}} \sum_{|\alpha|=m} |z^\alpha|^p \leq \lambda^p  \norm{Q}^p_{B_{Z}}  (|z_{1}|^p + \cdots + |z_{n}|^p)^m,$$ 
		which is less than or equals to $ \lambda^p  \norm{Q}^p_{B_{Z}} \left(\sup_{z \in B_{Z}}\norm{z}_{p}\right)^{pm}$.
		\comment{\begin{align} \label{e-1.9}
				\sum_{|\alpha|=m} \norm{U(a_{\alpha})(rz)^{\alpha}}^p \nonumber
				&= r^{pm} \sum_{|\alpha|=m} \norm{U(a_{\alpha})z^{\alpha}}^p\\ \nonumber
				& \leq r^{pm} \norm{U}^p  \norm{Q}^p_{B_{Z}} \sum_{|\alpha|=m} |z^\alpha|^p \\ \nonumber
				& \leq r^{pm} \norm{U}^p  \norm{Q}^p_{B_{Z}}  (|z_{1}|^p + \cdots + |z_{n}|^p)^m \\ \nonumber
				& = r^{pm} \norm{U}^p  \norm{Q}^p_{B_{Z}} \norm{z_{p}}^{pm}\\ 
				& \leq r^{pm} \norm{U}^p  \norm{Q}^p_{B_{Z}} \left(\sup_{z \in B_{Z}}\norm{z_{p}}\right)^{pm}\\ \nonumber
				& \leq r^{pm} \lambda^p  \norm{Q}^p_{B_{Z}} \left(\sup_{z \in B_{Z}}\norm{z_{p}}\right)^{pm},
			\end{align}
			which is less than or equal to $\lambda^p \norm{Q}^p_{B_{Z}}$ if
			$$
			r \leq \frac{1}{\sup_{z \in B_{Z}}\norm{z_{p}}}.
			$$}Therefore, 
		$K^m_{\lambda}(B_{Z}, p,U) \geq \frac{1}{\sup_{z \in B_{Z}}\norm{z}_{p}}$,
		 which together \eqref{e-1.8} and Lemma \ref{lem-3.3}(1) yields
		\begin{equation} \label{e-1.13}
			R_{\lambda}(B_{Z}, p,U) \geq \frac{1}{4 \lambda\,2^{\frac{1}{p}}}\left(\frac{\lambda^p - \norm{U}^p}{2\lambda^p - \norm{U}}\right)^{\frac{1}{p}} \frac{1}{\sup_{z \in B_{Z}}\norm{z}_{p}}.
		\end{equation}	
		On the other hand, for any $z \in B_{Z}$, we have $$\sum_{|\alpha|=m} \norm{U(a_{\alpha}) z^{\alpha}}^p \leq  \norm{U}^p  \norm{Q}^p_{B_{Z}} \left(\sup_{z \in B_{Z}}\norm{z}_{p}\right)^{pm}.
		$$ 
		Thus, $K^m(B_{Z}, p,U) \geq \frac{1}{\norm{U}^{\frac{1}{m}}\,\sup_{z \in B_{Z}}\norm{z}_{p}}.
		$
		Using this from Lemma \ref{lem-3.3}(2) and \eqref{e-1.8}, we deduce that 
		\begin{equation} \label{e-1.14}
			R_{\lambda}(B_{Z}, p,U) \geq \frac{1}{4 \lambda\,2^{\frac{1}{p}}}\left(\frac{\lambda^p - \norm{U}^p}{\lambda^p - \norm{U}^p +1}\right)^{\frac{1}{p}}\, \frac{1}{\norm{U}} \frac{1}{\sup_{z \in B_{Z}}\norm{z}_{p}}, \,\, \mbox{if} \,\, \norm{U} \geq \frac{1}{4 \lambda\,2^{\frac{1}{p}}}
		\end{equation}
		and
		\begin{equation} \label{e-1.15}
			R_{\lambda}(B_{Z}, p,U) \geq \left(\frac{\lambda^p - \norm{U}^p}{\lambda^p - \norm{U}^p +1}\right)^{\frac{1}{p}}\, \frac{1}{\sup_{z \in B_{Z}}\norm{z}_{p}}, \,\, \mbox{if} \,\, 0<\norm{U} <\frac{1}{4 \lambda\,2^{\frac{1}{p}}}.
		\end{equation}
		Finally, the conclusion follows from \eqref{e-1.13}, \eqref{e-1.14}, and \eqref{e-1.15}. This completes the proof.
	\end{proof}
	To simplify computations in the proofs of the remaining theorems, we may assume without loss of generality that all pluriharmonic functions $f \in \mathcal{PH}(B_{Z},X)$ satisfy $\norm{f}_{B_{Z},X} \leq 1$, where $X=\mathcal{B}(\mathcal{H})$. In other words, we restrict attention to the class $$\mathcal{BPH}(B_{Z}, X):=\{f: f \in \mathcal{PH}(B_{Z},X)\, \mbox{ with}\, \norm{f}_{B_{Z},X} \leq 1\}.$$
	 This is justified since any $f\in \mathcal{PH}(B_{Z},X)$ with $\norm{f}_{B_{Z},X} > 1$ can be rescaled by refining the coefficients. Indeed, if $f$ has the form \eqref{e-1.3-a} with $\norm{f}_{B_{Z},X} >1$, then defining $F$ by $F(z)=f(z)/\norm{f}_{B_{Z},X}$ yields $\norm{F}_{B_{Z},X} =1$ with coefficients $a_{\alpha}$, $b_{\alpha}$ replaced by $a_{\alpha}/\norm{f}_{B_{Z},X}$ and $b_{\alpha}/\norm{f}_{B_{Z},X}$ in \eqref{e-1.3-a}. In the same reasoning, we consider 
	 $$\mathcal{BPH}(^mB_{Z},X):=\{P: P \in \mathcal{PH}(^mB_{Z},X)\,\, \mbox{with}\,\, \norm{P}_{B_{Z},X} \leq 1\}.
	 $$
	\begin{pf} [{\bf Proof of Theorem \ref{thm-1.2}}]
		We want to first obtain a lower bound of $R^{m}_{\lambda}(B_{Z}, p,X)$ when $X$ is finite dimensional. For all $m\geq1$, let $P(z)=\sum_{|\alpha|=m} a_{\alpha}\, z^{\alpha} + \sum_{|\beta|=m} b^{*}_{\alpha}\, \bar{z}^{\alpha} \in \mathcal{BPH}(^mB_{Z},X)$. The proof is presented in three separate cases, depending on the values of $p$: $p=1$, $p\geq 2$, and $1<p<2$. Let $H$ and $G$ be the functions as defined in the proof of Theorem \ref{thm-1.1}. \\
		\underline{Case $p=1$:} 
		For each $\phi \in B_{X^{*}}$, we now consider the holomorphic polynomial $\tilde{H}:=\phi \circ H:\mathbb{C}^n \rightarrow \mathbb{C}$ defined by
		$\tilde{H}(z)=\phi(H(z))= \sum_{|\alpha|=m} \phi\left(\frac{a_{\alpha} + b_{\alpha}}{2}\right) z^{\alpha}$.
		Clearly, $|\tilde{H}(z)| \leq \norm{H(z)}$ for each $z$. Since $X$ is finite dimensional, the identity operator $I_{X}$ is $1$-summing. Then for all $z \in B_{Z}$, we have $$M_{H,1}:=\sum_{|\alpha|=m} \norm{\left(\frac{a_{\alpha} + b_{\alpha}}{2}\right)z^{\alpha}} \leq \pi_{1}(I_{X}) \sup_{\phi \in B_{X^{*}}} \left(\sum_{|\alpha|=m}\left|\phi\left(\frac{a_{\alpha} + b_{\alpha}}{2}\right) z^{\alpha}\right|\right),$$
		 which is again less than or equals to 
		 $$\pi_{1}(I_{X}) \sup_{\phi \in B_{X^{*}}} \sup_{z \in B_{Z}} \left|\sum_{|\alpha|=m}\left|\phi\left(\frac{a_{\alpha} + b_{\alpha}}{2}\right) \right|z^{\alpha}\right|.$$
		  This quantity by the definition of $\chi_{M}(\mathcal{P}(^m Z))$ is less than or equals to 
		  $$\pi_{1}(I_{X})\, \chi_{M}(\mathcal{P}(^m Z)) \, \sup_{\phi \in B_{X^{*}}} \sup_{z \in B_{Z}} \left|\sum_{|\alpha|=m}\phi\left(\frac{a_{\alpha} + b_{\alpha}}{2}\right) z^{\alpha}\right|.
		  $$ 
		  Thus, 
		  $M_{H,1} \leq \pi_{1}(I_{X}) \chi_{M}(\mathcal{P}(^m Z)) \sup_{\phi \in B_{X^{*}}}\sup_{z \in B_{Z}} |\tilde{H}(z)| \leq \pi_{1}(I_{X})\, \chi_{M}(\mathcal{P}(^mZ)) \,   \norm{H}_{B_{Z},X}.
		  $
		\comment{\begin{align} \label{e-1.19}
				\sum_{|\alpha|=m} \norm{\left(\frac{a_{\alpha} + b_{\alpha}}{2}\right)z^{\alpha}} \nonumber
				&\leq \pi_{1}(I_{X}) \sup_{\phi \in B_{X^{*}}} \left(\sum_{|\alpha|=m}\left|\phi\left(\frac{a_{\alpha} + b_{\alpha}}{2}\right) z^{\alpha}\right|\right) \\ \nonumber
				& \leq \pi_{1}(I_{X}) \sup_{\phi \in B_{X^{*}}} \sup_{z \in B_{Z}} \left(\sum_{|\alpha|=m}\left|\phi\left(\frac{a_{\alpha} + b_{\alpha}}{2}\right) z^{\alpha}\right|\right) \\ \nonumber
				& = \pi_{1}(I_{X}) \sup_{\phi \in B_{X^{*}}} \sup_{z \in B_{Z}} \left|\sum_{|\alpha|=m}\left|\phi\left(\frac{a_{\alpha} + b_{\alpha}}{2}\right) \right|z^{\alpha}\right| \\ \nonumber
				& \leq \pi_{1}(I_{X})\, \chi_{M}(\mathcal{P}(^m Z)) \, \sup_{\phi \in B_{X^{*}}} \sup_{z \in B_{Z}} \left|\sum_{|\alpha|=m}\phi\left(\frac{a_{\alpha} + b_{\alpha}}{2}\right) z^{\alpha}\right| \\ \nonumber
				& = \pi_{1}(I_{X})\, \chi_{M}(\mathcal{P}(^m Z)) \, \sup_{\phi \in B_{X^{*}}} \sup_{z \in B_{Z}} \left|\phi \left(\sum_{|\alpha|=m}\left(\frac{a_{\alpha} + b_{\alpha}}{2}\right) z^{\alpha}\right)\right| \\ 
				& = \pi_{1}(I_{X})\, \chi_{M}(\mathcal{P}(^mZ)) \,   \norm{H}_{B_{Z},X} 
			\end{align}
			where the fourth inequality follows from the definition of $\chi_{M}(\mathcal{P}(^m Z))$.} 
		By considering the holomorphic polynomial $\tilde{G}:=\phi \circ G:\mathbb{C}^n \rightarrow \mathbb{C}$, proceeding as above we obtain $M_{G,1} \leq \pi_{1}(I_{X})\, \chi_{M}(\mathcal{P}(^m Z)) \,\, \norm{G}_{B_{Z},X}$ for all $z \in B_{Z}$. By making use of these inequalities, and then from Lemma \ref{lem-3.1}, for all $z \in B_{Z}$, a simple computation shows that 
		$$\sum_{|\alpha|=m} \norm{a_{\alpha}z^{\alpha}} \leq \pi_{1}(I_{X})\, \chi_{M}(\mathcal{P}(^m Z)) \, \left(\norm{H}_{B_{Z},X}+\norm{G}_{B_{Z},X}\right)  \leq 4\, \pi_{1}(I_{X})\, \chi_{M}(\mathcal{P}(^m Z)).
		$$
		Similarly, for all $z \in B_{Z}$, we have $\sum_{|\alpha|=m} \norm{b_{\alpha}z^{\alpha}} \leq 4\, \pi_{1}(I_{X})\, \chi_{M}(\mathcal{P}(^m Z))$.
		The last two inequalities give 
		$$\sum_{|\alpha|=m} (\norm{a_{\alpha}}+\norm{b_{\alpha}})|z^{\alpha}| \leq 8\, \pi_{1}(I_{X})\, \chi_{M}(\mathcal{P}(^m Z)),
		$$
		which yields 
		$$R^{m}_{\lambda}(B_{Z}, p,X)  \geq E_{1}(X)\, \lambda^{\frac{1}{m}} \, \left(\chi_{M}(\mathcal{P}(^m Z)) \right)^{-\frac{1}{m}},$$
		where $E_{1}(X)$ is a constant depending on $X$. Plugging the estimates from \cite[(4.5), (4.6), pp. 187]{defant-2003} in the above inequality, we obtain 
		\begin{equation*}
			R^{m}_{\lambda}(B_{Z}, p,X) \geq E_{1}(X)\, \lambda^{\frac{1}{m}} \, \max \left\{\frac{1}{\sqrt{n}\, \norm{Id: \ell^n_{2} \rightarrow Z}}, \, \frac{1}{\left(\frac{m^m}{m!}\right)^{1/m}\, \norm{Id: Z \rightarrow \ell^n_{1}}}\right\},
		\end{equation*}
		and hence the desired lower bound of $R_{\lambda}(B_{Z}, p,X)$ follows from Lemma \ref{lem-3.3} (1) and using the fact $\sup _{m \in \mathbb{N}} (m^m/m!)^{1/m} =e$. This completes the proof for $p=1$.
		\\ [1mm]
		\underline{Case $p\geq 2$:} Let $p\geq 2$. Note that  $$M_{H,p}=\sum_{|\alpha|=m} \norm{\left(\frac{a_{\alpha} + b_{\alpha}}{2}\right)z^{\alpha}}^p \leq \sum_{|\alpha|=m} \norm{\left(\frac{a_{\alpha} + b_{\alpha}}{2}\right)z^{\alpha}}^2 \leq  \sum_{|\alpha|=m} \norm{\left(\frac{a_{\alpha} + b_{\alpha}}{2}\right)}^2
		$$
		for any $z \in B_{Z}$. Hence, by the fact $\norm{(a_{\alpha}+b_{\alpha})/2}=\sup_{\phi \in B_{X^{*}}}|\phi((a_{\alpha}+b_{\alpha})/2)|$ and the estimate in \cite[p. 187]{defant-2003}, we obtain 
		$M_{H,p}\leq  \norm{I:\ell^n _2 \rightarrow Z}^{2m} \norm{H}^2_{B_{Z},X}$.
		Similarly, we have $M_{G,p} \leq \norm{I:\ell^n _2 \rightarrow Z}^{2m} \norm{G}^2_{B_{Z},X}$.
		\comment{\begin{align} \label{e-1.24}
				\sum_{|\alpha|=m} \norm{\left(\frac{a_{\alpha} + b_{\alpha}}{2}\right)z^{\alpha}}^p \nonumber
				& \leq \sum_{|\alpha|=m} \norm{\left(\frac{a_{\alpha} + b_{\alpha}}{2}\right)z^{\alpha}}^2 \\ \nonumber
				& \leq  \sum_{|\alpha|=m} \norm{\left(\frac{a_{\alpha} + b_{\alpha}}{2}\right)}^2 \\ 
				& \leq  \norm{I:\ell^n _2 \rightarrow Z}^{2m} \norm{H}^2_{B_{Z},X}.
			\end{align}
			Similarly,
			\begin{equation} \label{e-1.25}
				\sum_{|\alpha|=m} \norm{\left(\frac{a_{\alpha} - b_{\alpha}}{2}\right)z^{\alpha}}^p \leq \norm{I:\ell^n _2 \rightarrow Z}^{2m} \norm{G}^2_{B_{Z},X}.
		\end{equation}}
		Then these two inequalities together with Lemma \ref{lem-3.1} yield 
		$$\left(\sum_{|\alpha|=m} \norm{a_{\alpha}z^{\alpha}}^p\right)^{\frac{1}{p}}  \leq 2 \norm{I:\ell^n _2 \rightarrow Z}^{\frac{2m}{p}}  2^{\frac{2}{p}}.$$
		This inequality remains valid if the left hand quantity replaced by $\left(\sum_{|\alpha|=m} \norm{b_{\alpha}z^{\alpha}}^p\right)^{1/p} $.
		A little calculation using these inequalities 
		shows 
		$$\sum_{|\alpha|=m} \left(\norm{a_{\alpha}}^p + \norm{b_{\alpha}}^p\right) |z^{\alpha}|^p \leq 2^{p+3}  \, \norm{I:\ell^n _2 \rightarrow Z}^{2m},$$
		and as a consequence, 
		$R^{m}_{\lambda}(B_{Z}, p,X) \geq 2^{-\frac{p+3}{pm}}\,  \lambda^{\frac{1}{m}}\, \norm{I:\ell^n _2 \rightarrow Z}^{-\frac{2}{p}}$.
		Finally, the desired lower bound of $R_{\lambda}(B_{Z}, p,X)$ follows from Lemma \ref{lem-3.3} (1).\\
		\underline{Case $1<p<2$:} For any $1<p<2$, using H\"{o}lder's inequality we obtain 
		\begin{equation*}
			M_{H,p} \leq  \left(\sum_{|\alpha|=m} \norm{\left(\frac{a_{\alpha} + b_{\alpha}}{2}\right)z^{\alpha}}\right)^{2-p} \left(\sum_{|\alpha|=m} \norm{\left(\frac{a_{\alpha} + b_{\alpha}}{2}\right)z^{\alpha}}^2\right)^{p-1}.
		\end{equation*}
		Lemma \ref{lem-3.1}, combined with the bounds of $M_{H,1}$ and $M_{H,2}$ from the above two cases, yields 
		$$M_{H,p}\leq 2^{p}\,(\pi_{1}(I_{X}))^{2-p}\, (\chi_{M}(\mathcal{P}(^mZ)))^{2-p} \norm{I:\ell^n _2 \rightarrow X}^{2(p-1) m}.$$
		\comment{\begin{align*}
				&\sum_{|\alpha|=m} \norm{\left(\frac{a_{\alpha} + b_{\alpha}}{2}\right)z^{\alpha}}^p \\ \nonumber
				&=\sum_{|\alpha|=m} \norm{\left(\frac{a_{\alpha} + b_{\alpha}}{2}\right)z^{\alpha}}^{p(1-\theta)}\, \norm{\left(\frac{a_{\alpha} + b_{\alpha}}{2}\right)z^{\alpha}}^{p\theta} \\ \nonumber
				& \leq \left(\sum_{|\alpha|=m} \norm{\left(\frac{a_{\alpha} + b_{\alpha}}{2}\right)z^{\alpha}}\right)^{p(1-\theta)} \left(\sum_{|\alpha|=m} \norm{\left(\frac{a_{\alpha} + b_{\alpha}}{2}\right)z^{\alpha}}^2\right)^{\frac{p\theta}{2}} \\ \nonumber
				& \leq (\pi_{1}(I_{X}))^{p(1-\theta)}\, (\chi_{M}(\mathcal{P}(^mZ)))^{p(1-\theta)} \norm{I:\ell^n _2 \rightarrow X}^{p \theta m} \norm{H}^{p}_{B_{Z},X} \\ \nonumber
				& \leq 4^{p}\,(\pi_{1}(I_{X}))^{p(1-\theta)}\, (\chi_{M}(\mathcal{P}(^mZ)))^{p(1-\theta)} \norm{I:\ell^n _2 \rightarrow X}^{p \theta m} \, \norm{I-\real(a_{0})}^{p}\, \norm{f}^p_{B_{Z},X}.
		\end{align*}}This last inequality remains valid if we replace the left-hand side by $M_{G,p}$ as well. Thus, we obtain $$\left(\sum_{|\alpha|=m} \norm{a_{\alpha}z^{\alpha}}^p\right)^{1/p} \leq 
		4\,(\pi_{1}(I_{X}))^{(2-p)/p}\, (\chi_{M}(\mathcal{P}(^mZ,Y)))^{(2-p)/p} \norm{I:\ell^n _2 \rightarrow X}^{ 2(p-1)/(pm)}.$$
		 This inequality remains valid if the left-hand side is replaced by $\left(\sum_{|\alpha|=m} \norm{b_{\alpha}z^{\alpha}}^p\right)^{1/p}$. Using these, a simple calculation shows that 
		$$\sum_{|\alpha|=m} \left(\norm{a_{\alpha}}^p + \norm{b_{\alpha}}^p\right) |z^{\alpha}|^p \leq 2^{2p+1}(\pi_{1}(I_{X}))^{2-p}\, (\chi_{M}(\mathcal{P}(^mZ,Y)))^{2-p} \norm{I:\ell^n _2 \rightarrow X}^{2(p-1) m},$$
		and thus, 
		$R^{m}_{\lambda}(B_{\ell^n _q}, p,X) \geq E_{2}(X)\, \lambda^{\frac{1}{m}} \,(\chi_{M}(\mathcal{P}(^mZ)))^{-\frac{(2-p)}{pm}} \norm{I:\ell^n _2 \rightarrow X}^{-\frac{2(p-1)}{p}}$, where $E_{2}(X)$ is a constant depending on $X$ and $p$. Plugging the estimates from \cite[(4.5), (4.6), pp. 187]{defant-2003} in the last inequality and letting $\theta = (2(p-1))/p$, we obtain 
		\begin{equation*}
			R^{m}_{\lambda}(B_{Z}, p,X) \geq E_{2}(X)\, \lambda^{\frac{1}{m}} \, \max \left\{\frac{1}{(\sqrt{n})^{1-\theta}\, \norm{Id: \ell^n_{2} \rightarrow Z}}, \, \frac{\norm{I:\ell^n _2 \rightarrow Z}^{-\theta}}{\left(\frac{m^m}{m!}\right)^{(1-\theta)/m}\, \norm{Id: Z \rightarrow \ell^n_{1}}^{1-\theta}}\right\},
		\end{equation*}
		and hence the desired lower bound of $R_{\lambda}(B_{Z}, p,X)$ follows from Lemma \ref{lem-3.3} (1) and using the fact $\sup _{m \in \mathbb{N}} (m^m/m!)^{1/m} =e$. We now want to find the desired upper bound of $R_{\lambda}(B_{Z}, p,X)$. In view of Lemma \ref{lem-3.5}, it is enough to know the upper bound of $R_{\lambda}(\Omega, p,X)$ for a complete Reinhardt domain $\Omega$. For $\Omega=B_{\ell^n_{q}}$, from \cite{A-H-P-Forum}, it is known that 
		\begin{equation*}
			K_{\lambda}(B_{\ell^n_{q}}, p,X) \leq d\, \lambda^{\frac{2}{\log \, n}}\, n^{1/q}\,\frac{\left(\log n\right)^{1-\left(1/\min\{q, 2\}\right)}}{n^{\frac{1}{p}+\frac{1}{\max\{2, q\}}-\frac{1}{2} }}=d\, \lambda^{\frac{2}{\log \, n}}\,n^{1-\frac{1}{p}}\, \left(\frac{\log \, n}{n}\right)^{1-\left(1/\min\{q, 2\}\right)}
		\end{equation*}
		for some constant $d>0$. Now the desired upper bound of $R_{\lambda}(B_{Z}, p,X)$ follows from the fact $R_{\lambda}(B_{\ell^n_{q}}, p,X) \leq K_{\lambda}(B_{\ell^n_{q}}, p,X)$ and from Lemma \ref{lem-3.5}. This completes the proof.
	\end{pf}
	\begin{pf} [{\bf Proof of Theorem \ref{thm-1.3-a}}]
		We first recall the following well-known observations: if $Y=(\mathbb{C}^n, ||.||)$ is a Banach space for which the canonical basis vectors $e_{k}$ form a normalized $1$-unconditional basis then  $(i)$ $\norm{Id: Y \rightarrow \ell_2^n}=\sup_{z \in B_{Y}}(\sum_{m=1}^{n}|z_{m}|^2)^{1/2}$, $(ii)$ $\norm{Id: Y \rightarrow \ell_1^n}=\sup_{z \in B_{Y}} \sum_{m=1}^{n}|z_{m}|=\norm{\sum_{m=1}^{n}e^{*}_{m}}_{Y^{*}}$ for all $z \in Y$. We now come to the proof of this theorem.\\
		\underline{Proof of part $(1)$:} Since $Z \subset \ell_2$, we may assume that the inclusion map $Z \subset \ell_2$ has norm $1$. The desired conclusion then follows from Theorem \ref{thm-1.2}, together with the preceding facts $(i)$ and $(ii)$.\\
		\underline{Proof of part $(2)$:} In view of \cite[Theorem 1.d.5]{Lindenstrauss-book}, we may assume that $M^{(2)}(X)=1$. Consequently, $\norm{Id: \ell_2^n \rightarrow Z_{n}}=1$. The proof now follows from Theorem \ref{thm-1.2}, the above facts $(i)$ and $(ii)$, and \cite[(5.3), p. 191]{defant-2003}.This completes the proof.	
	\end{pf}
	\begin{pf} [{\bf Proof of Theorem \ref{thm-1.3}}]
		We first want to find the lower bound. We split the proof into the following cases. Let $H$ and $G$ be the functions as defined in the proof of Theorem \ref{thm-1.1}. Let $P$ be as in Theorem \ref{thm-1.2}.\\
		\underline{Case $p > q$:} Since $p > q$, we have $M_{H,p}
		\leq M_{H,q}$ for $z \in B_{\ell^n _q}$. So, estimate in \cite[(3
		.7), p. 24]{defant-2006} give
		\begin{align*} 
			M_{H,p}
			\leq  \sum_{|\alpha|=m} \left(e^{\frac{m}{q}} \left(\frac{m!}{\alpha!}\right)^{\frac{1}{q}}\norm{H}_{B_{\ell^n _q},X}\right)^q  |z^\alpha|^q  
			= e^m \norm{H}^q_{B_{\ell^n _q},X} \big(\sum_{j=1}^{n}|z_{j}|^q\big)^m \leq e^m \norm{H}^q_{B_{\ell^n _q},X}
		\end{align*}
		for any $z \in B_{\ell^n _q}$.
		Similarly, $M_{G,p} \leq e^m  \norm{G}^q_{B_{\ell^n _q},X}$. By Lemma \ref{lem-3.1}, using the last two inequalities we deduce that $\left(\sum_{|\alpha|=m} \norm{a_{\alpha}z^{\alpha}}^p\right)^{1/p} \leq 2\, e^{\frac{m}{p}}\, 2^{\frac{q}{p}}$. This inequality remains valid if the left-hand side is replaced by $\left(\sum_{|\alpha|=m} \norm{b_{\alpha}z^{\alpha}}^p\right)^{1/p}$.
	By making use of these inequalities, we obtain 
	$\sum_{|\alpha|=m} \left(\norm{a_{\alpha}}^p + \norm{b_{\alpha}}^p\right) |z^{\alpha}|^p \leq 2^{p+1}\, e^m \, 2^q$,
	 which shows that 
	 $R^{m}(B_{\ell^n _q}, p,X) \geq 2^{-\left(1+\frac{1}{p}\right)\frac{1}{m}}\, e^{-\frac{1}{p}}\, 2^{-\frac{q}{pm}}$. Finally, in view of Lemma \ref{lem-3.3} (2), we deduce that
	\begin{equation*}
		R_{\lambda}(B_{\ell^n _q}, p,X) \geq \dfrac{\left(\lambda^p -1\right)^{\frac{1}{p}}}{\lambda}\, 2^{-\left(1+\frac{1}{p}\right)}\, e^{-\frac{1}{p}}\, 2^{-\frac{q}{p}}.
	\end{equation*}
	\underline{Case $p > Cot(X)$:} In view of Lemma \ref{lem-3.4}, to find a lower bound of $R_{\lambda}(B_{\ell^n _q}, p,X)$ it is enough to find that of $R_{\lambda}(\mathbb{D}^n, p,X)$. Let $f \in \mathcal{BPH}(\mathbb{D}^n,X)$ be of the form \eqref{e-1.3-a}. For any complex Banach space $X$ with a finite cotype $t$, by making use of \cite[Theorem 3.1]{carando-2020} to the function $H$ we have
	\begin{equation} \label{e-1.35}
		\left(\sum_{|\alpha|=m} \norm{\left(\frac{a_{\alpha}+b_{\alpha}}{2}\right)}^t\right)^{\frac{1}{t}} \leq c^m \, \left(\int\limits_{\mathbb{T}^n} \norm{H(z)}^t \, dz\right)^{\frac{1}{t}} \leq c^m \norm{H}_{\mathbb{D}^n,X}
	\end{equation}
	for some constant $c>0$. Now, for any $p > Cot(X)$, there exists $t \in [Cot(X), p)$ such that $X$ has cotype $t$, and hence $X$ has cotype $p$. Using \eqref{e-1.35}, we have $\sum_{|\alpha|=m} \norm{\left(a_{\alpha}+b_{\alpha}\right)/2}^p \leq c^{mp} \norm{H}^p_{\mathbb{D}^n,X}$.
	Similarly, by considering the function $G$, there exists constant $d>0$ such that $\sum_{|\alpha|=m} \norm{\left(a_{\alpha}-b_{\alpha}\right)/2}^p \leq d^{mp} \norm{G}^p_{\mathbb{D}^n,X}$.
	For any $p\geq 1$, using last two inequalities from Lemma \ref{lem-3.1} we deduce that $\sum_{|\alpha|=m} \norm{a_{\alpha}}^p \leq 2^p (c^m + d^m)^p \,\, \mbox{and} \,\, \sum_{|\alpha|=m} \norm{b_{\alpha}}^p \leq 2^p (c^m + d^m)^p$. By these inequalities, for $z \in \mathbb{D}^n$, we obtain $$\sum_{|\alpha|=m} \left(\norm{a_{\alpha}}^p + \norm{b_{\alpha}}^p\right) |z^{\alpha}|^p \leq \sum_{|\alpha|=m} \left(\norm{a_{\alpha}}^p + \norm{b_{\alpha}}^p\right) \leq 2^{p+1} \,(c^m + d^m)^p,$$
	which shows 
	$R^{m}(\mathbb{D}^n, p,X) \geq E_{5}\, 2^{-\left(1+\frac{1}{p}\right)\frac{1}{m}}$,
	 where $E_{5}>0$ is a constant independent of $n$. In view of Lemma \ref{lem-3.3} (2) and Lemma \ref{lem-3.4} we get
	\begin{equation*}
		R_{\lambda}(B_{\ell^n _q}, p,X) \geq R_{\lambda}(\mathbb{D}^n, p,X) \geq E_{5} \,2^{-\left(1+\frac{1}{p}\right)}\, \dfrac{\left(\lambda^p -1\right)^{\frac{1}{p}}}{\lambda}.
	\end{equation*}
	\underline{Case $p \leq \min \{Cot(X),q\}$:} Let $p \leq \min\{Cot(X),q\}$, then $p \leq t$ if $X$ has cotype $t$. Let $\beta:=\min\{t,q\}$. Let $f \in \mathcal{BPH}(B_{\ell^n _q},X)$ be of the form \eqref{e-1.3-a}. Then from the above two cases, for any $z \in B_{\ell^n _q}$, there exists constant $c>0$ such that 
	$\sum_{m=0}^{\infty}\sum_{|\alpha|=m} (\norm{a_{\alpha}}^\beta + \norm{b_{\alpha}}^\beta) | (\delta\,z)^{\alpha}|^\beta \leq \lambda^{\beta}$,
	where $\delta= c\, (\lambda^p - 1)^{1/p}/\lambda$. Using this inequality together with H\"{o}lder's inequality, and letting $r=\gamma \, \delta \in (0,1)$, we obtain $\sum_{m=0}^{\infty}\sum_{|\alpha|=m} (\norm{a_{\alpha}}^p + \norm{b_{\alpha}}^p) | (r\,z)^{\alpha}|^p$ is less than or equals to
	$$
	\sum_{m=0}^{\infty} \gamma^{pm} ((\sum_{|\alpha|=m}\norm{a_{\alpha}(\delta z)}^{\beta})^{p/\beta}+(\sum_{|\alpha|=m}\norm{b_{\alpha}(\delta z)}^{\beta})^{p/\beta})(\sum_{|\alpha|=m}1)^{1-p/\beta} \leq \lambda^p \sum_{m=0}^{\infty} \big(c_{1}\gamma  n^{\left(1/p-1/\beta\right)}\big)^{pm},$$ 
	for some constant $c_{1}>0$, independent of $n$. Now the last quantity is less than or equals to $\lambda^p$ if $\gamma= c_{2}\, n^{\left(1/p-1/\beta\right)}$, where $c_{2}$ is a constant, independent of $n$. Hence, we conclude that $$R_{\lambda}(B_{\ell^n _q}, p,X) \geq  E_{6} \,\, \dfrac{\left(\lambda^p -1\right)^{1/p}}{\lambda}\,\, n^{\left({\frac{1}{\beta}}-\frac{1}{p}\right)}$$
	for some constant $E_{6}>0$, independent of $n$.
	Now, we want to find the upper bound of $R_{\lambda}(B_{\ell^n _q}, p,X)$. To prove this, we shall make use of Theorem \ref{thm-A}. By considering $z=e_{j}$ in Theorem \ref{thm-A}, for each $n \in \mathbb{N}$ there exist $a_{1}, \ldots, a_{n} \in X$ such that $1/2 \leq \norm{a_{j}}$ for each $1\leq j \leq n$. Using this inequality and considering the $1$-homogeneous polynomial $F \in \mathcal{PH}(^1B_{Z},X)$ of the form $F(z)=\sum_{j=1}^{n} z_{j}a_{j}+ \sum_{j=1}^{n} a^{*}_{j}\bar{z_{j}}$, from Theorem \ref{thm-A} and the definition of $R^1_{\lambda}(B_{\ell^n _q}, p)$, a simple computation shows that
	\begin{align} \label{e-1.39}
		\frac{2\,n}{2^p\, n^{\frac{p}{q}}} \leq \sum_{j=1}^{n}(\norm{a_{j}}^p+\norm{a_{j}}^p) n^{-\frac{p}{q}} \leq \frac{\norm{F}_{B_{\ell^n_q},X}}{(R^1_{\lambda}(B_{\ell^n _q}, p,X))^p} \leq \frac{2}{(R^1_{\lambda}(B_{\ell^n _q}, p,X))^p}\, \sup_{z \in B_{\ell^n _q}} \norm{z}^p_{Cot(X)}.
	\end{align} 
	For $q \leq Cot(X)$, we have $\sup_{z \in B_{\ell^n _q}} \norm{z}_{Cot(X)} \leq 1$, which by \eqref{e-1.39} give 
	$R^1(B_{\ell^n _q}, p,X)\leq 2\, n^{1/q\, - \,1/p}.$ On the other hand, for $q>Cot(X)$, $\sup_{z \in B_{\ell^n _q}} \norm{z}_{Cot(X)} \leq n^{1/Cot(X)\, - \,1/p}$, and so from \eqref{e-1.39} we obtain $R^1(B_{\ell^n _q}, p) \leq 2\, n^{1/Cot(X)\, - \,1/p}$. It is known that $R_{\lambda}(B_{\ell^n _q}, p,X)\leq  R^1_{\lambda}(B_{\ell^n _q}, p,X)=\lambda\,R^1(B_{\ell^n _q}, p,X)$, and hence the upper bound of $R_{\lambda}(B_{\ell^n _q}, p,X)$ immediately follows. This completes the proof.
\end{pf}

\begin{proof}[{\bf Proof of Theorem \ref{thm-1.1-a}}]
	In the proof of Theorem \ref{thm-1.1}, we have shown that for $X=\mathcal{B}(\mathcal{H})$, $K^m_{\lambda}(B_{Z}, p,U) \geq \frac{1}{\sup_{z \in B_{Z}}\norm{z}_{p}}$ and $K^m(B_{Z}, p,U) \geq \frac{1}{\norm{U}^{\frac{1}{m}}\,\sup_{z \in B_{Z}}\norm{z}_{p}}$, and the same reasoning extends to any complex Banach space $X$. The desired bound of $K_{\lambda}(B_{Z}, p,U) $ then follows from Remark \ref{rem-4.1}, completing the proof.
\end{proof}

\noindent{\bf Statements and Declarations:}\\

\noindent{\bf Acknowledgment.} The author gratefully acknowledges the support received during the initial stages of this work while serving as a Senior Project Associate at IIT Bombay, through the \textit{Core Research Grant (CRG)}, sanctioned to Professor Sourav Pal by the Science and Engineering Research Board (SERB), Government of India. The author is currently supported by the J. C. Bose Fellowship of Professor B. V. Rajarama Bhat, ISI Bangalore, whose support is also sincerely acknowledged.\\

\noindent{\bf Conflict of interest.} The author declares that there is no conflict of interest regarding the publication of this paper.\\

\noindent{\bf Data availability statement.} Data sharing not applicable to this article as no datasets were generated or analysed during the current study.\\

\noindent{\bf Competing Interests.} The author declares none.


\end{document}